\documentclass[11pt]{article}

\usepackage{amsmath,amssymb,amsthm}
\usepackage{booktabs}
\usepackage{graphicx}
\usepackage{xcolor}
\usepackage{geometry}
\usepackage[colorlinks=true,linkcolor=blue!70!black,citecolor=green!50!black,urlcolor=blue!70!black,
            pdftitle={High-Precision Approximation of Riemann Zeros via the Truncated Weil Form},
            pdfauthor={Akiva Groskin},
            pdfsubject={Riemann Hypothesis; Connes-van Suijlekom truncation; Galerkin numerical analysis},
            pdfkeywords={Riemann zeta function, Weil quadratic form, Connes-van Suijlekom truncation, Galerkin method, Sobolev regularity, Dirichlet L-functions, arbitrary-precision computation, independent reproduction}]{hyperref}
\usepackage{url}
\makeatletter
\g@addto@macro\UrlBreaks{\do\/\do-\do\_\do.\do\:}
\makeatother

\theoremstyle{plain}
\newtheorem{theorem}{Theorem}[section]

\theoremstyle{definition}

\theoremstyle{remark}
\newtheorem{remark}[theorem]{Remark}
\newtheorem{heuristic}[theorem]{Heuristic}

\newcommand{\dps}{\mathrm{dps}}

\newcommand{\CvS}{Connes--van Suijlekom}
\newcommand{\CCM}{Connes--Consani--Moscovici}
\newcommand{\CC}{Connes--Consani}
\DeclareMathOperator{\keff}{k_{eff}}

\title{High-Precision Approximation of Riemann Zeros\\
  via the Truncated Weil Form}

\author{Akiva Groskin}

\date{May 2026}

\begin{document}
\maketitle

\noindent
\textbf{MSC\,2020:} 11M26 (primary), 65F15, 65N25 (secondary).

\medskip
\noindent
\textbf{Keywords:} Riemann zeta function, Weil quadratic form,
Connes--van Suijlekom truncation, Galerkin method,
Sobolev regularity, Dirichlet $L$-functions,
arbitrary-precision computation, independent reproduction.

\begin{abstract}
The \CvS\ truncated Weil quadratic form, indexed by a cutoff
parameter~$c$ that controls the primes $p \leq c$ entering the
operator, produces a ground state whose Fourier--Mellin zeros
provably lie on the critical line; whether they converge to the
Riemann zeros as $c\to\infty$ is open (Connes 2026; \CCM\ 2025).
We present, to our knowledge, the first independent public
implementation of the \CvS\ Galerkin matrix at sixteen cutoffs
($c=13$ through $67$, plus $c=100$).  Across the in-sample window
$c=13$ through $c=67$ at $N=100$, the first-zero absolute error
$|\gamma_1 - \gamma_1^{\mathrm{Riemann}}|$ shrinks monotonically
from $\sim\!2\times10^{-55}$ to $\sim\!1.5\times10^{-168}$ ---
a $113$-OOM convergence across fifteen cutoffs
(Table~\ref{tab:15pt-sweep}).  The smallest-positive even-sector
eigenvalue $\lambda_{\min}^{\mathrm{even}}$ separately reaches
${\sim}10^{-334}$ at $c=100$, $N=250$ ($275$-OOM span from $c=13$).

\emph{Out-of-sample test at $c=100$.}  On the four-point $N$-sweep
$N\in\{100,150,200,250\}$ at $\dps=500$, consecutive first-difference
ratios $0.837$ and $0.836$ match to two decimal places.
Aitken-$\Delta^2$ on the two overlapping triples yields
$\log_{10}|\lambda_\infty^{\mathrm{even}}|\approx -536.8$ and
$\approx -533.7$, approaching the Connes 2026 \S6.4 heuristic
prediction ($\approx -530.4$) monotonically with $N$ ($6.4$ and
$3.3$~OOM gaps out of $|x_\infty|\sim 530$).  The same eigenvector
recovers $\gamma_1,\ldots,\gamma_{10}$ to $307$--$329$ matching
digits at $N=250$, $\dps=500$.  Under the unitary equivalence with
\CCM\ Lemma~5.1, this is the deepest such Galerkin-truncation
recovery in the public \CvS/\CCM\ literature.  The
\emph{smallest-positive} even-sector eigenvalue is the genuine
smallest eigenvalue: a small block of negative-sign eigenvalues seen
at the finite archimedean cutoff $T=800$ is an artifact of that
cutoff and is absent once $T$ is increased (continuum positivity
of $QW_\lambda$ is RH-equivalent and is not assumed at
$\lambda=\sqrt{100}$).  The fit
$|\log_{10}\lambda_{\min}|\approx 13.24\,c^{0.634}$ on $c\leq 67$
at $N=100$ is shown to be a finite-$N$ rate, falsified at
$c=100, N=200$ by $49$~OOM in the direction of faster decay.

\emph{Structural observations} include approximate eigenvector
$c$-invariance (overlap $\geq 0.9498$ on all $105$ cutoff pairs
despite eigenvalues differing by $113$~OOM), multi-zero
convergence universality (all ten detectable zeros within $3.8\%$
of each other), an empirical Galerkin-convergence exponent
$s(c)\approx 55\log c - 128$, un-rescaled Galerkin bulk-spectrum
Poisson statistics ($\beta<0.05$; this is a structural diagnostic
of the truncated operator, not a test of Montgomery's conjecture,
which applies to locally-rescaled zero spacings), and tight bulk
invariants
$\log|\det Q_c|\approx-65.6\,c+542$ ($R^2=0.997$).

We make no claim of proof; the contribution is reproducible
numerical data and its careful interpretation under the
existing CvS/CCM framework.

All code, data, and ancillary files are publicly available.
\end{abstract}

\medskip
\noindent\textbf{Note added in revision (2026 June).} The negative-sign
eigenvalue blocks reported at $c=100$ (Section~\ref{sec:c100-doublets})
and for $L(s,\chi_3)$ at $c=23,29$ (Section~\ref{sec:chi3-odd}) are
artifacts of the finite archimedean integration cutoff $T$: they are
stable in working precision but vanish once $T$ is increased, so
cutoff-free the relevant even sectors are non-negative and the
smallest-positive branch is the genuine smallest eigenvalue.  No
quantitative result changes; the $\gamma_k$ extractions, the Aitken
extrapolation, and all convergence data are unaffected.  We thank A.~Connes, whose questions about the $c=100$ spectrum prompted
this investigation; the cutoff sensitivity was independently identified
by B.~W.~A.~Silva, consistent with the naturally even, positive ground
state reported by R.~Andrews (full attributions in the Acknowledgments).

\medskip
\noindent\textbf{Note added (2026 August).} Three further corrections, none
of which changes a measured value.  (i)~The Paley--Wiener mechanism
proposed in Section~\ref{sec:sobolev-pw} is \emph{withdrawn}.  The test
pre-registered for it in Section~\ref{sec:future} has been carried out: at
$c=23$ on the $N\in\{40,60,80\}$ grid at $\dps=150$, the empirical exponent
is $s=46.140$, $46.031$ and $45.934$ at $T=400$, $800$ and $1600$, so a
fourfold increase in~$T$ moves~$s$ by $0.206$, where the proposed relation
$s=\sigma_{\mathrm{eff}}T/2\pi$ requires proportionality in~$T$ at fixed~$c$.
By the criterion stated in Section~\ref{sec:future} the mechanism is
therefore incorrect, and the constant $A=0.432$ has no meaning as derived.
Relatedly, Table~\ref{tab:sobolev-scaling} is measured at $T=400$, inherited
from the per-cutoff $N$-convergence ladders
(Table~\ref{tab:N-conv}), and its caption now says so;
Section~\ref{sec:sobolev-pw} had substituted the $T=800$ of the fifteen-cutoff
sweep.  The empirical scaling $s(c)\approx 55\log c-128$ is a fit to the
measured exponents rather than a consequence of the withdrawn mechanism, and
is unaffected, as are the measured exponents themselves.  The three-cutoff
test was first carried out by M.~Osman and has been reproduced independently
here.  (ii)~Under the matching-digit definition given in the caption of
Table~\ref{tab:c100-gamma}, the $c=67$, $N=100$, $\dps=200$ datum, whose
error is $1.478\times10^{-168}$, gives $167$ matching digits, not the $168$
previously printed in Section~\ref{sec:c100-gamma}; the error value itself is reported
correctly throughout.  (iii)~Two summary ranges in
Section~\ref{sec:c100-gamma} were misreported: the increments are $93$--$117$
matching digits for $\dps$ $500\to1000$ at $N=150$, and $181$--$203$ for
$N$ $150\to250$ at $\dps=500$, not the $95$--$115$ and $179$--$201$ previously
printed.  The underlying per-$\gamma_k$ counts, of which these ranges are
differences, are unchanged and are reported correctly throughout.

\section{Introduction}\label{sec:intro}

In February 2026, Connes published ``Past, Present and a Letter Through
Time'' \cite{Connes2026}, revisiting the spectral interpretation of the
Riemann zeros within the framework of noncommutative geometry.
Section~6 of that paper presents what is, in effect, an open question:
for each value of a cutoff parameter~$c$ (controlling the primes
$p\leq c$ entering the operator; $c$ itself need not be prime), the
truncated Weil quadratic form on
$L^2([0,\log c])$ admits a unique even-sector ground state whose
Fourier--Mellin zeros provably lie on the critical line.  Do these zeros
converge to the actual Riemann zeros as $c\to\infty$?

The criticality statement is a theorem.  Connes--van Suijlekom
\cite{CvS2025} establish (Theorem~6.1) that the zeros of
$\hat{\eta}_c$ are real for every finite~$c$, so the open question is
specifically about convergence, not about criticality.  If the
convergence holds, RH follows by Hurwitz's theorem applied to the
Mellin transforms --- but the convergence itself is unproven.

Connes reports numerical data for the first fifty zeros at $c=13$,
with errors ranging from approximately $2.6\times10^{-55}$ (first zero)
to approximately $10^{-3}$ (fiftieth zero); the first-zero error
relative to $\gamma_1=14.134725\ldots$ is the datum reproduced in this
work.  \CCM\ \cite{CCM2025} report the related number
$2.44\times10^{-55}$ at the same cutoff, with $N=120$ and 200-digit
precision.  As of April 2026, no
independent public reproduction of these numbers existed, and no public
code implementing the \CvS\ Galerkin matrix had been located in a
literature sweep encompassing arXiv, GitHub, MathOverflow, and seminar
recordings.  Both statements have since been overtaken: the code accompanying
this paper is public, and five parties have reported reproductions, three of
them from implementations independent of that code
(Section~\ref{sec:independent-repro}).

\subsection*{Contributions}
This paper makes five primary contributions:
\begin{enumerate}
\item \textbf{Independent reproduction and extension.}  We provide,
  to our knowledge, the first public implementation of the
  \CvS\ Proposition~4.1 Galerkin matrix (Python/mpmath); reproduce
  the $c=13$ first-zero error within a factor of~1.3 of the value
  reported by Connes 2026; cross-validate against \CCM\ at $c=14$;
  and extend the measurements to fifteen cutoffs ($c=13$
  through~$67$) spanning $113$~orders of magnitude in
  $|\gamma_1\text{ error}|$, plus a sixteenth cutoff at $c=100$.

\item \textbf{Out-of-sample empirical test of Connes 2026 \S6.4 at
  $c=100$.}  At $c=100$ we measure
  $\log_{10}|\lambda_{\min}^{\mathrm{even}}|$ on a four-point $N$-sweep
  $N\in\{100,150,200,250\}$ at $\dps=500$ (with a $\dps=1000$ retest at
  $N=150$), with values $-190.92, -247.19, -294.31, -333.68$.  The
  consecutive first-difference ratios $0.8373$ and $0.8355$ match to
  two decimal places, evidence for a local geometric model of the
  convergence sequence.  Aitken-$\Delta^2$ applied to the two
  consecutive overlapping triples gives
  $\log_{10}|\lambda_\infty^{\mathrm{even}}(c=100)|\approx -536.76$
  (from $\{100,150,200\}$) and $\approx -533.70$ (from
  $\{150,200,250\}$); the Connes 2026 \S6.4 heuristic continuum
  prediction is $\approx -530.38$, giving gaps of $6.39$ and $3.32$
  orders of magnitude respectively, out-of-sample (the in-sample fit
  window was $c\leq67$).  Section~\ref{sec:c100-verification} presents
  the computation.

\item \textbf{Recognition that the numerical pipeline measures the
  \CCM\ operator's spectrum.}  Under the unitary equivalence
  established in \S\ref{sec:spectral-triple} with \CCM\ (2025)
  Lemma~5.1, our $F_{\mathrm{even}}$ test function is, up to a
  positive scaling constant, the Fourier--Mellin transform
  $\widehat{\xi}_N(z)$ appearing in \CCM\ Theorem~1.1(iii); we recognize
  --- we do not claim to discover --- that our extracted zeros are,
  modulo the hypothesis status at $c=100$
  (Remark~\ref{rem:ccm-hyp-c100}), eigenvalues of the \CCM\ rank-one
  operator $D_{\log}^{(\lambda,N)}$ at $\lambda=\sqrt{c}$.  This makes
  our fifteen-cutoff sweep plus the $c=100$ datum, to our knowledge,
  the first independent multi-$c$ survey of that spectrum at
  trigonometric-Galerkin matching-digit precision, beyond \CCM's own
  reported values at $c\in\{9,12,13,14\}$; the deepest $\gamma_1$
  extraction at $c=100$ reaches $329$~matching digits at $N=250$,
  $\dps=500$.

\item \textbf{Structural characterization of the finite-$N$ operator.}
  We identify quantitative properties that constrain theoretical work:
  (a)~the Sobolev regularity scales as $s(c)\approx55\cdot\log c-128$
  at six cutoffs;
  (b)~all ten detectable zeros converge at rates within $3.8\%$ of
  each other;
  (c)~the eigenvector is approximately $c$-invariant
  ($>0.950$ overlap across all $105$~pairs);
  (d)~the normalized zero-crossing shape varies by only $3.3\%$;
  (e)~interval extension, not prime identity, drives convergence
  ($r=-0.96$);
  (f)~the spectral trace decays as $c^{-0.116}$ while the Frobenius
  norm is constant;
  (g)~the bulk spectrum is Poisson (no level repulsion);
  (h)~nine convergence models fail
  (Section~\ref{sec:extension});
  (i)~the empirical $|\log_{10}\lambda|\approx 13.24\,c^{0.634}$ fit
  is a finite-$N=100$ rate, not the continuum asymptote, falsified at
  $c=100, N=200$ by $49$~OOM and corroborated by a $c=67, N=150$
  rerun showing a $46$~OOM drop versus the same-$c$, $N=100$ datum;
  (j)~a block of negative-sign eigenvalues at $c=100$ is an artifact
  of the finite archimedean cutoff $T=800$: it is stable in working
  precision ($\dps\in\{500,1000\}$) but not in $T$, vanishing once $T$
  is increased, so cutoff-free the $c=100$ even sector carries no
  negative eigenvalues (Section~\ref{sec:c100-doublets}).

\item \textbf{Partial extension to Dirichlet $L$-functions.}  A port
  to $L(s,\chi_3)$ converges at $c=13,17,37$ ($16$--$28$~digits); an
  apparent positivity breakdown at $c=23,29$ is an artifact of the
  finite archimedean cutoff $T=400$, absent once $T$ is increased
  (Section~\ref{sec:chi3-odd}).
\end{enumerate}

\subsection*{Explicit non-claims}
We do not claim to validate Connes' convergence conjecture.  We do not
claim a proof of the Riemann Hypothesis.  We do not claim that the
Aitken extrapolation at $c=100$ proves a convergence rate.  The
agreement at $c=13$ is a consistency check against published values,
not a comparison against an independent prediction.  The agreement at
$c=100$ between our Aitken extrapolation and the Connes 2026 \S6.4
prediction is an empirical consistency test, not a theorem; we are
explicit that Aitken acceleration carries assumptions about the local
form of the convergence sequence, and the Connes 2026 \S6.4 prediction
is itself a heuristic statement based on the agreement (their Figure~1)
between $1-\chi_2(\lambda)$ and $\varepsilon(\lambda)$ on
$\lambda\leq 14$.  We report the agreement; we do not promote it.

\section{Mathematical Setup}\label{sec:setup}

\subsection{The Weil quadratic form and its truncation}

We follow the notation of \CvS\ \cite{CvS2025} and \CC\ \cite{CC2020}.
For a cutoff parameter~$c$ (controlling the primes $p\leq c$ in the
prime-sum piece below; $c$ itself need not be prime, as in our
$c=14$ and $c=100$ cells), set $L=\log c$.  The quadratic form acts on
$L^2([0,L])$ and decomposes as
\begin{equation}\label{eq:D-decomp}
  D = D_\infty + D_{\mathrm{pole}} + D_{\mathrm{prime}},
\end{equation}
where $D_\infty$ is the archimedean contribution,
$D_{\mathrm{pole}}$ is the pole contribution from the zeta function,
and $D_{\mathrm{prime}}$ sums over primes $p\leq c$.

The archimedean Mellin multiplier is
\begin{equation}\label{eq:h-plus}
  h_+(\tau) = -\log\pi
    + \operatorname{Re}\Psi\!\Bigl(\tfrac{1}{4}+\tfrac{i\tau}{2}\Bigr),
\end{equation}
where $\Psi$ denotes the digamma function.  This follows from \CC\
\cite{CC2020} (Equation~153) and has been independently fact-checked
against the derivation chain in \CC\ Appendix~B.

\subsection{The \CvS\ Proposition~4.1 Galerkin matrix}

\CvS\ Proposition~4.1 defines a function
$\psi\colon\mathbb{R}\to\mathbb{R}$ via
\begin{equation}\label{eq:psi-def}
  \psi(x) = \frac{1}{\pi}\int_0^L
    \sin\!\bigl(2\pi x(1-y/L)\bigr)\,D(y)\,dy,
\end{equation}
and the Galerkin matrix entries
\begin{equation}\label{eq:q-mn}
  q_{m,n} = \frac{\psi(m)-\psi(n)}{m-n},
  \qquad
  q_{n,n} = \psi'(n).
\end{equation}
The basis is the complex Fourier system
$\{U_n(x)=L^{-1/2}\exp(2\pi inx/L)\}_{n\in\mathbb{Z}}$ per
\CvS\ Proposition~4.1 / eq.~(5); the sin functions appearing in the
kernel of $\psi(x)$ above are the integrand, not the basis.  The
even-sector projection under the involution $y\mapsto L-y$ is
described in Section~\ref{sec:eigensolver}; only that sector is
diagonalized.

\subsection{Unitary equivalence with the \CCM\ matrix}

\CCM\ \cite{CCM2025} (Lemma~5.1) define a matrix $\tau_{i,j}$ in
terms of a function~$b_n$:
\begin{equation}\label{eq:ccm-matrix}
  b_n = -\frac{1}{\pi}\int_0^L \sin(2\pi ny/L)\,\mathcal{D}(y)\,dy,
  \qquad
  \tau_{i,j} = \frac{b_i - b_j}{i-j}.
\end{equation}
The key identity connecting the two formulations is, for integer~$n$:
\begin{equation}\label{eq:sin-identity}
  \sin\!\bigl(2\pi n(1-y/L)\bigr) = -\sin(2\pi ny/L),
\end{equation}
which follows from the $2\pi n$-periodicity of sine.  Substituting
into~\eqref{eq:psi-def} at integer argument:
\begin{equation*}
  \psi(n)
  \;=\; \frac{1}{\pi}\int_0^L \!-\!\sin(2\pi ny/L)\,D(y)\,dy
  \;=\; -\frac{1}{\pi}\int_0^L \sin(2\pi ny/L)\,D(y)\,dy
  \;=\; b_n.
\end{equation*}
Therefore $b_n = \psi(n)$ and $\tau_{i,j} = (b_i - b_j)/(i-j)
= (\psi(i) - \psi(j))/(i-j) = q_{i,j}$ as matrices.  The two
matrices coincide identically (no sign flip): the eigenvalues and
eigenvectors of $\tau$ and $q$ agree, and the present work and the
\CCM\ numerics are computing the same operator via two equivalent
matrix presentations.

\subsection{Spectral-triple interpretation of \texorpdfstring{$F_{\mathrm{even}}$}{F\_even}}
\label{sec:spectral-triple}

A direct consequence of the unitary equivalence in \S2.3 connects the
$F_{\mathrm{even}}$ test function used throughout this work to the
rank-one spectral-triple construction of \CCM\ \cite{CCM2025}.  In
their Section~5, \CCM\ define a self-adjoint rank-one perturbation of
the scaling operator
\begin{equation}\label{eq:Dlog-perturbed}
  D_{\log}^{(\lambda,N)}
  \;:=\; D_{\log}^{(\lambda)}
  \;-\; \bigl|D_{\log}^{(\lambda)}\,\xi_N\bigr\rangle
        \bigl\langle\delta_N\bigr|,
\end{equation}
acting on $L^2([\lambda^{-1},\lambda], du/u)$ as a rank-one
perturbation of $D_{\log}^{(\lambda)}$, supported on the
$(2N{+}1)$-dimensional subspace $E_N$, with $\lambda^2=c$.
\CCM\ Theorem~1.1(iii) establishes that the spectrum of
$D_{\log}^{(\lambda,N)}$ coincides with the zeros of
\begin{equation}\label{eq:xi-hat}
  \widehat{\xi}_N(z)
  \;=\; \int \xi_N(u)\,u^{-iz}\,\frac{du}{u},
\end{equation}
the Fourier--Mellin transform of the ground-state eigenvector
$\xi_N$ of the truncated Weil quadratic form, taken with \CCM's
duality convention $\langle \mathbb{R}^*_+ \mid \mathbb{R}\rangle$.
Positivity of the continuum quadratic form $QW_\lambda$ for all
$\lambda > 1$ is equivalent to the Riemann Hypothesis
(Connes 2026~\cite{Connes2026} \S4.1, citing Weil) and is proved only
for small values of $\lambda$ (Yoshida, see~\cite{Connes2026} and
references therein).  We do not assume continuum positivity at
$\lambda = \sqrt{100}$.  Throughout this paper $\xi_N$ denotes the
\emph{empirically distinguished smallest-positive eigenvector} of the
finite-$N$ even-sector matrix $Q_c^{\mathrm{even}}(c,N)$ --- a
selection criterion, not a theorem-derived ground-state
identification.\footnote{At $c\leq67$ with $N=100$ the matrix is
observed positive on every cell, so the smallest eigenvalue coincides
with the smallest-positive one and our $\xi_N$ coincides with the
\CCM\ Theorem~1.1 ground state; the same holds cutoff-free at $c=100$ once the
finite-$T$ artifact of Remark~\ref{rem:ccm-hyp-c100} is removed, so
$\xi_N$ refers throughout to this smallest-positive ground state.}

\begin{remark}[Hypothesis status at $c=100$]
\label{rem:ccm-hyp-c100}
\CCM\ Theorem~1.1 takes its starting point as $\epsilon_N$, the
\emph{smallest eigenvalue} of $QW_\lambda^N$, assumed simple and
with even eigenvector.  At the cutoffs $c \leq 67$ studied in
Section~\ref{sec:extension}, the finite-$N$ Galerkin matrix
$Q_c^{\mathrm{even}}$ is observed positive on every cell at $N=100$,
its smallest eigenvalue is numerically simple and coincides with
its smallest-positive eigenvalue, and the corresponding eigenvector
projects predominantly onto the even sector --- so the
finite-dimensional data are consistent with the \CCM\ Theorem~1.1
ground-state hypotheses (as numerical conditions, not as
analytically established hypotheses of the continuum theorem).  At $c=100$ with $N\geq 100$
(Section~\ref{sec:c100-verification}), the raw finite-$N$ spectrum
computed at the finite archimedean cutoff $T=800$ carries a small
block of negative-sign eigenvalues (Section~\ref{sec:c100-doublets}).
These are an artifact of that cutoff: they are absent once $T$ is
increased, so cutoff-free the $c=100$ even sector is non-negative and
its smallest eigenvalue coincides with the smallest-\emph{positive}
branch we report, whose eigenvector gives the zero-extraction of
Table~\ref{tab:c100-gamma}.  The $c=100$ data are therefore consistent
with the \CCM\ Theorem~1.1 ground-state hypotheses on the same footing
as the $c\leq67$ cells.  Empirically, the positive branch is well-behaved: its
eigenvector extracts $\gamma_1,\ldots,\gamma_{10}$ to
$307$--$329$ matching digits at $c=100$, $N=250$.  The $c=67$ corroborative datum
(Section~\ref{sec:c100-reframe}) at $N=150$ remains in the
matrix-positive regime, providing a same-$N$ control where the
corresponding finite-dimensional positivity, simplicity, and
even-sector numerical conditions remain satisfied.
\end{remark}

The $F_{\mathrm{even}}$ test function defined for zero-extraction
(\S\ref{sec:eigensolver}) is, up to a positive scaling constant and
the change of variable $u = e^{x}$, the Fourier--Mellin transform of
the same eigenvector $\xi_N$ that appears in
equation~\eqref{eq:xi-hat}; the zero sets of the two transforms agree
on the real line by reality of $\xi_N$.  We identify the two via this
construction and adopt the convention $\lambda^2 = c$ from \CCM\
throughout.  Under this identification, and modulo the
hypothesis-status caveat of Remark~\ref{rem:ccm-hyp-c100} at
$c=100$, every $\gamma_k$ extracted in this paper is, equivalently,
a zero of $\widehat{\xi}_N$, and via \CCM\ Theorem~1.1(iii) an
eigenvalue of $D_{\log}^{(\lambda,N)}$ at $(\lambda{=}\sqrt{c},\, N)$.
The multi-cutoff data presented in Sections~\ref{sec:extension},
\ref{sec:c100-verification}, and~\ref{sec:structural} thereby provide
measurements of this spectrum at the cutoffs
$c\in\{13,\dots,67,100\}$, beyond \CCM's reported values at
$c\in\{9,12,13,14\}$ with $N=120$.

\subsection{Decomposition of the kernel}

The integral defining $\psi(x)$ decomposes into three pieces
corresponding to the three terms of~$D$:

\begin{enumerate}
\item \textbf{Archimedean piece} ($D_\infty$): involves the
  Mellin-multiplier integral $\int_0^T h_+(\tau)\,K(x,\tau)\,d\tau$
  where $K$ is a trigonometric kernel.  This piece admits a closed-form
  expression for the $-\log\pi$ component and requires numerical
  quadrature for the digamma component.

\item \textbf{Pole piece} ($D_{\mathrm{pole}}$): a rank-one correction
  arising from the simple pole of $\zeta(s)$ at $s=1$.

\item \textbf{Prime piece} ($D_{\mathrm{prime}}$): a finite sum over
  primes $p\leq c$ involving $\log p$ and the von~Mangoldt-type kernel.
\end{enumerate}

The reader needs to see these explicit pieces because the precision
characteristics of each --- the archimedean piece dominates the
quadrature cost, the pole piece is exact, and the prime piece involves a
finite sum --- determine the implementation strategy described in
Section~\ref{sec:implementation}.

\section{Numerical Implementation}\label{sec:implementation}

\subsection{Software stack and precision}

The implementation uses Python~3 with mpmath for arbitrary-precision
arithmetic.  The 15-cutoff $c\leq 67$ sweep of
Section~\ref{sec:extension} uses $\dps=80$ at the early cutoffs,
$\dps=150$ across the bulk, and $\dps=200$ at $c\geq 41$ to stay clear
of the backward-error floor (see Section~\ref{sec:precision}).  The
$c=100$ measurements of Section~\ref{sec:c100-verification} use
$\dps=500$ with a $\dps=1000$ retest at $N=150$.  Parallelism is
achieved via Python's \texttt{multiprocessing} module with 12-way
process parallelism.

The production entry points for the 15-cutoff sweep and the $c=100$
experiment are ancillary files \path{A_extended_c_sweep.py} and
\path{c100_experiment_optimized.py}, both included with this
submission.  The python-flint \texttt{acb.digamma} primitive provides
an $\sim 11\times$ speedup over pure mpmath in the archimedean
quadrature.  All ancillary scripts run against the public
\texttt{connes-cvs} package on PyPI (version 0.2.2 at the time of the
original submission); see the package version note in
Section~\ref{sec:code}.

\subsection{Stable kernel evaluation}

The trigonometric kernel arising in the quadrature of the archimedean
piece involves expressions of the form
\begin{equation}\label{eq:unstable-kernel}
  \frac{e^{i\beta L}-1}{i\beta}
\end{equation}
which suffer catastrophic cancellation when $|\beta|$ is small.  We
replace this with the numerically stable identity
\begin{equation}\label{eq:stable-kernel}
  \frac{e^{i\beta L}-1}{i\beta}
  = \frac{\sin(\beta L)}{\beta}
  + \frac{2i\sin^2(\beta L/2)}{\beta},
\end{equation}
together with a Taylor-series fallback (5~terms at $\dps=80$, 7~terms
at $\dps=150$) when $|\beta|$ falls below a threshold~$\varepsilon$.
The stable kernel has been verified against analytic values to
$10^{-136}$ accuracy at $\beta=10^{-55}$.

\subsection{Quadrature strategy}\label{sec:quadrature}

The Mellin-multiplier integral over $\tau\in[0,T]$ with $T=800$ is
evaluated using mpmath's adaptive tanh--sinh quadrature
(\texttt{mp.quad}, Takahasi--Mori~1974).  The $-\log\pi$ component of $h_+$ admits a
closed-form antiderivative, eliminating quadrature error for that piece.
The digamma component
$\operatorname{Re}\Psi(1/4+i\tau/2)$ is the dominant cost: each
evaluation requires one mpmath digamma call at approximately 4.5\,ms/call
at $\dps=150$.

The cutoff $T=800$ was chosen based on $T$-convergence testing
(see Section~\ref{sec:T-convergence}).
The $-\log\pi$ component is handled analytically (exact), and the
remaining digamma component decays as $O(\log\tau/\tau)$ against the
oscillatory kernel; the Riemann--Lebesgue lemma ensures convergence.
The agreement with Connes' and \CCM's independently computed $c=13$
datum (within a factor of~1.3) provides empirical evidence that
$T=800$ is sufficient, as any $O(1)$ truncation error in the matrix
would produce $O(1)$ eigenvalue errors incompatible with the observed
agreement.

\subsection{Even-sector projection and eigensolver}\label{sec:eigensolver}

The implementation uses the full Fourier basis
$\{e^{2\pi iky/L}\}_{k=-N}^{N}$ (dimension $2N+1=201$ at $N=100$),
which decomposes into even (cosine-type) and odd (sine-type) sectors
under the involution $y\mapsto L-y$.  The even sector comprises the
constant mode ($k=0$) and the $N$~cosine combinations
$(\mathbf{e}_k+\mathbf{e}_{-k})/\sqrt{2}$ for $k=1,\dots,N$,
yielding dimension $N+1=101$ at $N=100$.
Only the even sector is diagonalized.  The matrix parity error is zero by construction (the
off-diagonal blocks between even and odd sectors vanish identically for
the \CvS\ quadratic form).

The eigenvalue decomposition uses \texttt{mpmath.eigsy} in symmetric
mode.  The smallest \emph{positive} eigenvalue
$\lambda_{\min}^{\mathrm{even}}$ of the even-sector matrix is the
quantity of primary interest: it is the empirically distinguished
positive branch whose eigenvector gives the zero-extraction reported
throughout (see Remark~\ref{rem:ccm-hyp-c100} for the hypothesis-status
discussion).  At all fifteen cutoffs in Section~\ref{sec:extension} at
$N=100$ the matrix is observed positive on every cell, so smallest-overall
and smallest-positive coincide.  At $c=100$, $N\geq 100$
(Section~\ref{sec:c100-verification}) the raw matrix computed at the
finite cutoff $T=800$ carries a small block of negative-sign
eigenvalues (Section~\ref{sec:c100-doublets}) that are absent at larger
$T$, so cutoff-free the smallest-positive selection coincides with the
matrix-level smallest eigenpair.

\subsection{Root extraction}

The first Riemann zero $\gamma_1=14.134725141734693790\ldots$ is
located as the first zero of the Mellin transform of the ground-state
eigenvector.  Root extraction uses \texttt{mp.findroot} with the
Anderson solver and explicit tolerance $\mathrm{tol}=10^{-140}$ for
the general $c \leq 67$ sweep; the $c=100$ headline extraction at
$N=250$, $\dps=500$ uses the retightened tolerance $10^{-380}$
reported in Table~\ref{tab:c100-gamma}, sufficient to certify the
$307$--$329$ matching-digit counts.  The absolute error is
$|\gamma_1^{\mathrm{exact}}-\gamma_1^{\mathrm{computed}}|$.

\section{Results: Reproduction at \texorpdfstring{$c=13$}{c=13}}
\label{sec:reproduction}

\subsection{Primary measurement}

At $c=13$, $N=100$, $T=800$, $\dps=200$, the even-sector minimum
eigenvalue and first-zero error are:
\begin{equation}\label{eq:c13-result}
  \lambda_{\min}^{\mathrm{even}} = 2.865\times10^{-59},
  \qquad
  |\gamma_1\text{ error}| = 2.005\times10^{-55}.
\end{equation}
The measurement is stable across iterations: the Iteration~4
($N=100$, stable kernel) and Iteration~5 ($c=13$ cell) results agree
to leading digits.  The $\dps=150$ and $\dps=200$ retests reproduce
this cell, indicating that $\dps=80$ is fully sufficient at $c=13$.

\subsection{Comparison with published values}

\begin{table}[ht]
\centering
\caption{Comparison of $|\gamma_1\text{ error}|$ at $c=13$ across
  sources.}
\label{tab:c13-comparison}
\begin{tabular}{@{}lllll@{}}
\toprule
Source & $|\gamma_1\text{ error}|$ & Basis & $N$ & Precision \\
\midrule
Connes 2026, Section~6 \cite{Connes2026}
  & $\approx2.6\times10^{-55}$ & trigonometric & --- & --- \\
CCM 2511.22755, Intro \cite{CCM2025}
  & $\approx2.5\times10^{-55}$ & trigonometric & --- & --- \\
CCM 2511.22755, \S6 ($\lambda=\sqrt{13}$) \cite{CCM2025}
  & $2.44\times10^{-55}$ & trigonometric & 120 & 200 digits \\
This work
  & $2.005\times10^{-55}$ & trigonometric & 100 & 150 digits \\
\bottomrule
\end{tabular}
\end{table}

All four measurements cluster within a factor of~1.3.  The three
values from the Connes/CCM group cluster around $2.5\times10^{-55}$;
this work's value differs by at most that factor of~1.3 from any of
them, and by approximately~1.2 from the most precisely-reported
\CCM\ datum $2.44\times10^{-55}$.

\subsection{Discussion of the factor-of-1.3 discrepancy in the size of the
\texorpdfstring{$|\gamma_1\text{ error}|$}{|gamma\_1 error|}}
\label{sec:c13-1p3-discrepancy}

To be explicit: the factor of 1.3 refers to a discrepancy in the
\emph{size of the error term} $|\gamma_1 - \gamma_1^{\mathrm{Riemann}}|$,
not to a discrepancy in $\gamma_1$ itself nor in $\lambda_{\min}$.  All
four cited measurements report the same Riemann zero $\gamma_1$ to far
more than the leading digits of the error magnitude; the factor of 1.3
quantifies how much the implementations differ in how small the
residual error becomes after the Galerkin truncation.

The discrepancy is consistent with the differences in implementation:
\begin{itemize}
\item \textbf{Basis consistency with CCM.}  The \CCM\ Galerkin
  computation and this work both use the trigonometric basis: \CCM\
  Lemma~5.1 defines the matrix entries via the kernel
  $\sin(2\pi ny/L)$, and Connes 2026 \S6 describes the truncation as
  acting on the ``trigonometric orthonormal basis.''  The factor-1.3
  spread in the size of $|\gamma_1\text{ error}|$ across the four
  measurements is therefore not attributable to a basis-choice
  difference; it reflects differences
  in $N$, precision, integration cutoff~$T$, and normalization
  conventions, as detailed below.  A prolate-basis implementation of
  the CvS Galerkin would be a genuinely distinct path in the literature
  --- potentially valuable as an independent cross-check --- but is not
  yet reported by any group.
\item \textbf{Truncation parameter~$N$.}  \CCM\ uses $N=120$; this
  work uses $N=100$.  Differences in $N$ at the same basis choice
  produce different Galerkin truncation errors; at $c=13$, the data
  suggest that $N=100$ at 150-digit precision is sufficient to reach
  the precision floor.
\item \textbf{Precision.}  \CCM\ uses 200-digit precision; this work
  uses 200~digits for the final sweep.  At $c=13$, the agreement
  to all reported digits between $\dps=80$ and $\dps=150$ shows
  that the result is not precision-limited.
\item \textbf{Quadrature cutoff~$T$ and integration method.}  This work
  uses $T=800$; the \CCM\ quadrature details are not fully specified.
  Differences in quadrature strategy could contribute.
\end{itemize}

We emphasize that all four numbers are measurements of the same
quantity --- the distance between the first Fourier--Mellin zero of
$\eta_{13}$ and $\gamma_1$ --- via the same operator (up to the sign
equivalence of Section~\ref{sec:setup}) computed in the same
trigonometric basis, with different $N$, precision, and integration
parameters.  The factor-of-1.3 spread in the size of
$|\gamma_1\text{ error}|$ is the current state of the art for this
measurement.

\subsection{Convergence in \texorpdfstring{$N$}{N} at
  \texorpdfstring{$c=13$}{c=13}}

The project's iteration history provides a convergence sequence in~$N$:

\begin{table}[ht]
\centering
\caption{$|\gamma_1\text{ error}|$ at $c=13$, $T=400$, $\dps=80$,
  varying~$N$.}
\label{tab:c13-N-convergence}
\begin{tabular}{@{}lll@{}}
\toprule
$N$ & $|\gamma_1\text{ error}|$
    & $\lambda_{\min}^{\mathrm{even}}$ \\
\midrule
30  & $3.355\times10^{-44}$ & --- \\
60  & $4.22\times10^{-55}$  & --- \\
100 & $1.455\times10^{-55}$ & $2.077\times10^{-59}$ \\
\bottomrule
\end{tabular}
\end{table}
\noindent
(This measurement used $T=400$ quadrature; the $T=800$ value in Table~\ref{tab:15pt-sweep} ($2.865\times10^{-59}$) supersedes it.)

The $N=30\to60$ step yields an 11-order-of-magnitude improvement; the
$N=60\to100$ step yields a factor of approximately~2.9.  This rapid
initial convergence followed by saturation is consistent with the
eigenvector's effective support being concentrated in approximately 45
modes (see Section~\ref{sec:structural}): once $N$ exceeds
$k_{\mathrm{eff}}\approx45$, additional basis functions contribute
diminishing improvements.

\section{Extension to Fifteen Cutoffs}\label{sec:extension}

\subsection{Primary measurements at \texorpdfstring{$T=800$, $\dps=200$}{T=800, dps=200}}

The initial three-point sweep ($c\in\{13,17,19\}$) was performed at
$\dps=80$ (Iteration~5) and retested at $\dps=150$ (Iteration~6).  An
$11\times$ speed optimization (via python-flint's
\texttt{acb.digamma}) enabled the full sweep at $\dps=150$ in
Iteration~7.  Cutoffs $c\geq41$ were computed at $\dps=200$ to avoid
precision-floor contamination.  The final dataset uses $T=800$ at all
fifteen cutoffs.

\begin{table}[ht]
\centering
\caption{Complete 15-point $c$-sweep at $N=100$, $T=800$,
  $\dps=150$ ($c\leq37$) or $\dps=200$ ($c\geq41$).
  Rows $c\geq17$ have no prior published data we have located;
  the $c=14$ datum is, to our knowledge, the first independent
  measurement (CCM report a value at $\lambda=\sqrt{14}$).}
\label{tab:15pt-sweep}
\begin{tabular}{@{}rllllr@{}}
\toprule
$c$ & $L=\log c$ & $\lambda_{\min}^{\mathrm{even}}$
    & $|\gamma_1\text{ error}|$
    & $|\gamma_1\text{ err}|/\lambda_{\min}$
    & $\log_{10}|\gamma_1\text{ err}|$ \\
\midrule
13 & 2.565 & $2.865\times10^{-59}$
   & $2.005\times10^{-55}$  & 6999  & $-54.70$ \\
14 & 2.639 & $4.835\times10^{-65}$
   & $3.541\times10^{-61}$  & 7324  & $-60.45$ \\
17 & 2.833 & $2.030\times10^{-80}$
   & $1.634\times10^{-76}$  & 8047  & $-75.79$ \\
19 & 2.944 & $1.265\times10^{-90}$
   & $1.070\times10^{-86}$  & 8457  & $-85.97$ \\
23 & 3.135 & $5.959\times10^{-107}$
   & $5.520\times10^{-103}$ & 9262  & $-102.26$ \\
29 & 3.367 & $4.366\times10^{-124}$
   & $4.587\times10^{-120}$ & 10507 & $-119.34$ \\
31 & 3.434 & $1.045\times10^{-128}$
   & $1.141\times10^{-124}$ & 10919 & $-123.94$ \\
37 & 3.611 & $4.670\times10^{-140}$
   & $5.686\times10^{-136}$ & 12177 & $-135.25$ \\
41 & 3.714 & $2.122\times10^{-146}$
   & $2.760\times10^{-142}$ & 13004 & $-141.56$ \\
43 & 3.761 & $2.519\times10^{-149}$
   & $3.379\times10^{-145}$ & 13412 & $-144.47$ \\
47 & 3.850 & $2.994\times10^{-154}$
   & $4.270\times10^{-150}$ & 14260 & $-149.37$ \\
53 & 3.970 & $9.615\times10^{-161}$
   & $1.493\times10^{-156}$ & 15529 & $-155.83$ \\
59 & 4.078 & $2.328\times10^{-166}$
   & $3.911\times10^{-162}$ & 16800 & $-161.41$ \\
61 & 4.111 & $5.063\times10^{-168}$
   & $8.722\times10^{-164}$ & 17226 & $-163.06$ \\
67 & 4.205 & $7.993\times10^{-173}$
   & $1.478\times10^{-168}$ & 18489 & $-167.83$ \\
\bottomrule
\end{tabular}
\end{table}

\begin{figure}[!htbp]
\centering
\includegraphics[width=0.85\textwidth]{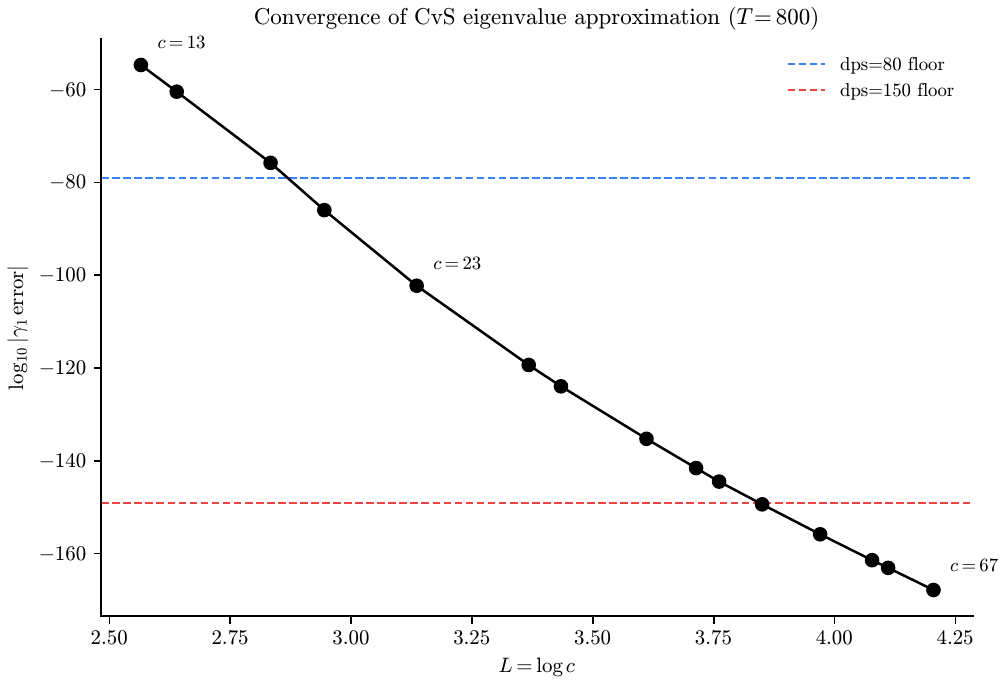}
\caption{First-zero absolute error $|\gamma_1\text{ error}|$ across
  fifteen cutoffs.  The data spans 113~orders of magnitude.
  Dashed lines indicate the backward-error floors at $\dps=80$ and
  $\dps=150$.}
\label{fig:convergence}
\end{figure}

The data spans 113~orders of magnitude in $|\gamma_1\text{ error}|$
from $c=13$ ($10^{-55}$) to $c=67$ ($10^{-168}$).  The convergence
is monotone (every increase in~$c$ produces a decrease in
$|\gamma_1\text{ error}|$) but the step sizes are non-uniform ---
see Section~\ref{sec:models} and Figure~\ref{fig:convergence}.

\begin{remark}[Precision-floor warning for $c=43$]\label{rem:c43}
At $\dps=150$, $\log_{10}|\gamma_1\text{ error}|=-144.63$ is only
approximately 1~order of magnitude above the $\dps=150$ backward-error
floor $\varepsilon\cdot\|Q\|_2\approx6\times10^{-150}$
($\log_{10}\approx-149.2$).  The $\dps=200$ spot-check
(Section~\ref{sec:precision}) corrected $\lambda_{\min}$ from
$1.836\times10^{-149}$ to $2.519\times10^{-149}$ and the
$|\gamma_1\text{ err}|/\lambda_{\min}$ ratio from 12900 to 13412, indicating that
the $\dps=150$ datum was mildly contaminated.  The $c=41$ datum at
$\log_{10}|\gamma_1\text{ error}|=-141.56$ has approximately 7~orders of
magnitude of headroom and is clean.
\end{remark}

\subsubsection{Cross-validation with \CCM\ at \texorpdfstring{$c=14$}{c=14}}

\CCM\ \cite{CCM2025} (Section~6) report
$|\gamma_1\text{ error}|\approx1.07\times10^{-60}$ at
$\lambda=\sqrt{14}$ (equivalently $c=14$), with $N=120$ and
200-digit precision.  This work measures
$3.541\times10^{-61}$ at the same cutoff --- a factor of
approximately~3 smaller.  This is the first row-for-row external
cross-validation beyond $c=13$.

The discrepancy is wider than the factor-of-1.3 observed at $c=13$.
The source of this discrepancy is not fully understood; plausible
contributions include $N$, precision, normalization conventions,
and integration parameters.

\subsection{\texorpdfstring{$T$}{T}-convergence at \texorpdfstring{$c=13$}{c=13}}
\label{sec:T-convergence}

To verify that the quadrature cutoff $T$ is sufficient, we compare
results at $c=13$ across two values of~$T$:

\begin{table}[ht]
\centering
\caption{$T$-convergence at $c=13$, $N=100$, $\dps=200$.}
\label{tab:T-convergence}
\begin{tabular}{@{}rlll@{}}
\toprule
$T$ & $\lambda_{\min}^{\mathrm{even}}$
    & $|\gamma_1\text{ error}|$
    & $\log_{10}|\gamma_1\text{ err}|$ \\
\midrule
400  & $2.077\times10^{-59}$ & $1.455\times10^{-55}$ & $-54.84$ \\
800  & $2.865\times10^{-59}$ & $2.005\times10^{-55}$ & $-54.70$ \\
\bottomrule
\end{tabular}
\end{table}

The shift from $T=400$ to $T=800$ is less than 0.14~units
in $\log_{10}|\gamma_1\text{ error}|$, indicating that the
quadrature is well-converged.  Both values remain within the
factor-of-1.3 envelope of the Connes/CCM measurements, indicating
that the $T=800$ result is not dominated by quadrature truncation
error.

\subsection{Pre-registered smooth-convergence model fits}
\label{sec:models}

Eight two-parameter smooth-convergence models were tested.
M1--M3 were pre-registered and fit to the initial 3-point $\dps=80$ data
(Iteration~5); M4--M5 were added for the 10-point sweep; M6--M8,
using prime-counting independent variables, were added for the full
15-point dataset.  All eight were fit to the 15-point
$T=800$ data.  The pre-registered acceptance criterion was
max residual $\leq0.5$ (in $\log_{10}$ units), corresponding to a
factor of $\sim\!3.2$ in the original error magnitude --- approximately
$0.4\%$ of the 113-OOM dynamic range.  All conclusions are robust to
relaxing the threshold to~2.0, at which point only M5 would pass.

\begin{align}
\text{M1 (linear in $L$):}&\quad
  \log_{10}|\mathrm{err}| = a\cdot L + b, \label{eq:M1} \\
\text{M2 (quadratic in $L$):}&\quad
  \log_{10}|\mathrm{err}| = a\cdot L^2 + b, \label{eq:M2} \\
\text{M3 ($L\log L$):}&\quad
  \log_{10}|\mathrm{err}| = a\cdot L\log L + b, \label{eq:M3} \\
\text{M4 (mixed linear $+$ log):}&\quad
  \log_{10}|\mathrm{err}| = a\cdot L + b\cdot\log L, \label{eq:M4} \\
\text{M5 (exponential):}&\quad
  \log_{10}|\mathrm{err}| = a\cdot\exp(b\cdot L), \label{eq:M5} \\
\text{M6 (linear in $\pi(c)$):}&\quad
  \log_{10}|\mathrm{err}| = a\cdot\pi(c) + b, \label{eq:M6} \\
\text{M7 (linear in $\vartheta(c)$):}&\quad
  \log_{10}|\mathrm{err}| = a\cdot\vartheta(c) + b, \label{eq:M7} \\
\text{M8 (prime powers):}&\quad
  \log_{10}|\mathrm{err}| = a\cdot\#\mathrm{pp}(c) + b. \label{eq:M8}
\end{align}

\begin{table}[ht]
\centering
\caption{Two-parameter smooth-convergence model fits to the 15-point
  data.}
\label{tab:model-fits}
\begin{tabular}{@{}llcc@{}}
\toprule
Model & Form & Max residual & Pass $\leq0.5$? \\
\midrule
M1 & $a\cdot L+b$            & 2.66  & FAIL \\
M2 & $a\cdot L^2+b$          & 4.71  & FAIL \\
M3 & $a\cdot L\log L+b$      & 3.63  & FAIL \\
M4 & $a\cdot L+b\cdot\log L$ & 2.66  & FAIL \\
M5 & $a\cdot\exp(b\cdot L)$  & 1.77  & FAIL \\
M6 & $a\cdot\pi(c)+b$        & 10.28 & FAIL \\
M7 & $a\cdot\vartheta(c)+b$  & 12.09 & FAIL \\
M8 & $a\cdot\#\mathrm{pp}(c)+b$ & 10.20 & FAIL \\
\bottomrule
\end{tabular}
\end{table}

All eight models are rejected.  M5 (exponential-in-$L$) is the best
among the $L$-based models at 1.77, but still exceeds the threshold
by a factor of~3.5.  M6--M8, which use prime-counting independent
variables ($\pi(c)$, Chebyshev~$\vartheta(c)$, and the number of
prime powers $\leq c$), perform substantially worse (max residuals
10--12), indicating that the convergence is not simply proportional to
any standard prime-counting function.

We emphasize that this failure does \emph{not} bear on the question of
whether the zeros converge as $c\to\infty$ --- it bears only on
whether the rate of convergence, if it exists, can be captured by a
two-parameter smooth model fit to fifteen points in the interval
$c\in[13,67]$.

\subsection{Structural observations}\label{sec:structural-obs}

Two structural features of the computation are worth recording.

\paragraph{Effective support invariance.}
Define $k_{\mathrm{eff}}(\epsilon)$ as the smallest~$k$ such that the
cumulative squared eigenvector weight
$S_k=\sum_{j=1}^k|v_j|^2$ satisfies $1-S_k<\epsilon$.  At
$\epsilon=10^{-50}$:

\begin{table}[ht]
\centering
\caption{Effective support $k_{\mathrm{eff}}(10^{-50})$ across
  cutoffs.}
\label{tab:keff}
\begin{tabular}{@{}llc@{}}
\toprule
$c$ & $\dps$ & $k_{\mathrm{eff}}(10^{-50})$ \\
\midrule
13 & 80/150 & 45 \\
17 & 150    & 45 \\
19 & 150    & 44 \\
\bottomrule
\end{tabular}
\end{table}

The ground-state eigenvector's effective support is approximately
44--47 modes regardless of cutoff.  This is approximately
$c$-invariant: the ``physical mode'' occupies a fixed number of basis
functions, and the remaining modes (approximately $k>75$ at $N=100$)
form a boundary-layer residual.

\paragraph{Eigenvector near-invariance across cutoffs.}
The overlap parameter $q=|\langle\eta_{c_1}|\eta_{c_2}\rangle|$ between
ground-state eigenvectors at different cutoffs satisfies $q>0.950$ for
all 105~measured cutoff pairs (see Section~\ref{sec:structural} for
the full analysis), despite eigenvalues differing by up to
113~orders of magnitude.  Three representative pairs from the initial
sweep:

\begin{table}[ht]
\centering
\caption{Eigenvector overlaps across cutoff pairs.}
\label{tab:overlaps}
\begin{tabular}{@{}lc@{}}
\toprule
Pair & $|\langle\eta_{c_1}|\eta_{c_2}\rangle|$ \\
\midrule
$(c=13,\;c=17)$ & 0.9977 \\
$(c=13,\;c=19)$ & 0.9956 \\
$(c=17,\;c=19)$ & 0.9997 \\
\bottomrule
\end{tabular}
\end{table}

What changes between cutoffs is predominantly the eigenvalue
$\lambda_{\min}$ (which shifts by up to 113~orders of magnitude), while the
eigenvector structure remains more than 95\% the same
(the minimum pairwise overlap across all 105~pairs is 0.950).

\subsection{Convergence in \texorpdfstring{$N$}{N} at fixed \texorpdfstring{$c=23$}{c=23}}
\label{sec:N-conv}

To characterize the Galerkin truncation error independently of the
prime-cutoff effect, we fix $c=23$ and vary the basis size
$N\in\{40,60,80,100,120,140\}$ at $\dps=150$.

\begin{table}[ht]
\centering
\caption{$N$-convergence at $c=23$, $T=400$, $\dps=150$.}
\label{tab:N-conv}
\begin{tabular}{@{}rllcr@{}}
\toprule
$N$ & $\lambda_{\min}^{\mathrm{even}}$
    & $|\gamma_1\text{ error}|$
    & $\keff(10^{-50})$
    & $\log_{10}|\gamma_1\text{ err}|$ \\
\midrule
 40 & $6.628\times10^{-72}$  & $1.033\times10^{-67}$  & 37 & $-67.0$ \\
 60 & $2.761\times10^{-88}$  & $3.044\times10^{-84}$  & 43 & $-83.5$ \\
 80 & $1.143\times10^{-99}$  & $1.113\times10^{-95}$  & 45 & $-95.0$ \\
100 & $4.182\times10^{-107}$ & $3.876\times10^{-103}$ & 46 & $-102.4$ \\
120 & $3.137\times10^{-111}$ & $2.854\times10^{-107}$ & 46 & $-106.5$ \\
140 & $9.036\times10^{-113}$ & $8.179\times10^{-109}$ & 46 & $-108.1$ \\
\bottomrule
\end{tabular}
\end{table}

The convergence in~$N$ is rapid for $N\leq80$ (+16.5 and
+11.4~orders of magnitude per 20-step increase) but saturates
sharply: the $N=100\to140$ range gains only 5.7~OOM total.
The effective support $\keff$ stabilizes at~46 from $N=100$ onward,
indicating that the ground-state eigenvector occupies approximately
46~even-sector basis functions regardless of the total basis size.
At $N=140$, the marginal improvement is only 1.5~OOM per 20-step ---
the Galerkin truncation has essentially converged, and the remaining
error of $8.18\times10^{-109}$ is controlled by the finite cutoff
$c=23$, not the finite basis.

\subsection{Log-periodic blind test}
\label{sec:logperiodic}

An empirical log-periodic model
\begin{equation}
\log_{10}|\gamma_1\text{ error}| = -66.95\cdot L + 110.98
  + 5.50\sin(3.08\cdot L + 6.57)
\end{equation}
was fit to the 10-point training data ($c=13$ through~$43$), yielding a
maximum training residual of~0.74 --- the first model to beat the
best two-parameter result of~1.77.  Pre-registered blind predictions were made for $c=47$
and~$c=53$ before computation.

\begin{table}[ht]
\centering
\caption{Log-periodic model blind test.}
\label{tab:logperiodic}
\begin{tabular}{@{}lrrrl@{}}
\toprule
$c$ & Predicted & Actual & $|\text{residual}|$ & Verdict \\
\midrule
47 & $-150.9$ & $-149.4$ & 1.5 & pass ($<2.0$) \\
53 & $-150.3$ & $-155.8$ & 5.5 & \textbf{fail} ($>2.0$) \\
\bottomrule
\end{tabular}
\end{table}

The model predicted that $c=53$ would yield \emph{less} precision
than $c=47$ (due to the oscillation swinging positive).  In reality,
$c=53$ showed $6.4$~OOM more precision.  The 5-parameter model was
overfitting the non-uniform step sizes in the training range; the
anti-monotone prediction was the tell.  This clean falsification
strengthens the negative result: no parametric model tested ---
including models with up to 5~free parameters --- captures the
$c$-dependence of $|\gamma_1\text{ error}|$ over the measured range.

\subsection{Blind prediction test at \texorpdfstring{$c=67$}{c=67}}
\label{sec:blind-c67}

Before computing the $c=67$ cell, an ensemble prediction (weighted
average of four independent prediction methods) was pre-registered:
\begin{equation}\label{eq:blind-pred}
  \log_{10}|\gamma_1\text{ error}|^{\mathrm{pred}} = -169.4 \pm 2.0.
\end{equation}
The prediction was sealed via SHA-256 hash before computation (see
ancillary file \path{preregistered_predictions_2026-04-12.md}).

The actual result was $\log_{10}|\gamma_1\text{ error}| = -167.83$.
The residual $|{-167.83} - ({-169.4})| = 1.57$ falls within the
pre-registered $\pm 2.0$ threshold: the blind prediction passes.

This is evidence that the convergence behavior is
sufficiently regular to permit out-of-sample predictions at 1.6-OOM
accuracy.  This is \emph{not} a proof of convergence; it is an
empirical observation that the convergence behavior is smooth and
predictable within the measured range.

\section{Out-of-Sample Empirical Test at \texorpdfstring{$c=100$}{c=100}}
\label{sec:c100-verification}

Operating at $c$ well above the precision floor enables a separate test
of the convergence behavior outside the in-sample window of
Section~\ref{sec:extension}.  The fifteen-cutoff sweep
($c\leq67$, $N=100$) established the empirical pattern
$\log_{10}|\lambda_{\min}^{\mathrm{even}}(c)| \approx -13.24\,c^{0.634}$
from a least-squares fit.  Two natural questions arise.  First, does
that fit, taken as an asymptotic, extrapolate correctly to larger $c$?
Second, is it consistent with the heuristic continuum asymptotic of
Connes 2026 \cite{Connes2026} \S6.4?  This section addresses both
through a measurement at $c=100$ with $N$-sweep at
$N\in\{100,150,200,250\}$ at $\dps=500$ and a precision retest at
$N=150$, $\dps=1000$.

\subsection{Connes 2026 Section~6.4 heuristic asymptotic}
\label{sec:c100-prediction}

Connes 2026 \cite{Connes2026} \S6.4 gives a heuristic continuum
asymptotic for the angular function $1-\chi_2(\lambda)$, motivated by
the agreement of $1-\chi_2(\lambda)$ with the smallest eigenvalue
$\varepsilon(\lambda)$ of the operator $A_\lambda$ for
$\lambda\leq14$ (their Figure~1).  The expression is
\begin{equation}\label{eq:connes-asymptotic}
  1-\chi_2(\lambda)
  \;\sim\; \frac{2^{14}}{3}\,\sqrt{2}\,\pi^{5}\;e^{-4\pi e^L + 9L/2},
\end{equation}
where, as defined in \cite{Connes2026} \S6.4, $L = 2\log\lambda$ is
the length of the support interval $[\lambda^{-1},\lambda]$ of the
admissible test functions for the truncated Weil form $QW_\lambda$.
With $\lambda^2=c$, this is $L = \log c$, so $e^L = c$.

Evaluating equation~\eqref{eq:connes-asymptotic} at $c=100$
($L=\log 100 \approx 4.6052$, $e^L=100$) in base-10 logarithm:
\begin{equation}\label{eq:connes-100-prediction}
\begin{aligned}
  \log_{10}\!\bigl(1-\chi_2(c=100)\bigr)
  &\approx
  \log_{10}\!\frac{2^{14}\,\sqrt{2}\,\pi^{5}}{3}
    - \frac{4\pi\cdot 100}{\ln10}
    + \frac{9\log 100}{2\ln10} \\
  &\approx 6.37 - 545.75 + 9.00 \;\approx\; -530.38.
\end{aligned}
\end{equation}
The prediction, treating the agreement of $\varepsilon(\lambda)$ with
$1-\chi_2(\lambda)$ on $\lambda\leq14$ as continuing to $c=100$, is
$\log_{10}|\varepsilon(c=100)|\approx -530.38$.

\subsection{Methodology}
\label{sec:c100-method}

We computed the Galerkin matrix $Q_c^{\mathrm{even}}$ at $c=100$,
$T=800$, working precision $\dps=500$, with python-flint
$\mathtt{flint\_prec}=2000$, using the v0.2.0 fused-kernel math core
from the public \texttt{connes-cvs} PyPI package, for
$N\in\{100,150,200,250\}$.  A precision retest at $\dps=1000$ with
$\mathtt{flint\_prec}=4000$ was run at $N=150$.

The retest provides three independent stability checks:
(i)~agreement of $\lambda_{\min}^{\mathrm{even}}$ between $\dps=500$
and $\dps=1000$ to $25$ leading significant digits;
(ii)~agreement of the absolute log-magnitudes of the five smallest
negative eigenvalues across $\dps=500$ and $\dps=1000$ to two decimal
places (these negatives are stable in $\dps$ but, as
Section~\ref{sec:c100-doublets} establishes, are a finite-$T$ cutoff
artifact);
(iii)~$\gamma_k$ extraction at $\dps=1000$ providing $\sim 219$--$242$
matching digits, doubling the depth at $\dps=500$.

Wall-clock at $c=100,N=150,\dps=500$ was 21.0~minutes on twelve workers
of an Apple~M2~Max; at $\dps=1000$, 111~minutes.

\subsection{Measurement: \texorpdfstring{$N$}{N}-sweep at \texorpdfstring{$c=100$}{c=100}}
\label{sec:c100-measurement}

\begin{table}[ht]
\centering
\caption{$N$-sweep at $c=100$, $T=800$.  Working precision $\dps=500$
for $N\in\{100,150,200,250\}$; retest at $\dps=1000$ for $N=150$.  The
\#~negative column counts finite-$T$ ($T=800$) artifact eigenvalues
(Section~\ref{sec:c100-doublets}), absent at larger $T$.}
\label{tab:c100-N-sweep}
\begin{tabular}{@{}rrcccr@{}}
\toprule
$N$ & $\dps$ &
$\lambda_{\min}^{\mathrm{even}}$ (smallest-positive) &
$\log_{10}|\lambda_{\min}|$ &
$\Delta$ in $\log_{10}$ &
\#~negative \\
\midrule
100 & 500  & $1.22 \times 10^{-191}$ & $-190.92$ & ---      & 3 \\
150 & 500  & $6.42 \times 10^{-248}$ & $-247.19$ & $-56.28$ & 5 \\
200 & 500  & $4.87 \times 10^{-295}$ & $-294.31$ & $-47.12$ & 8 \\
250 & 500  & $2.08 \times 10^{-334}$ & $-333.68$ & $-39.37$ & 11 \\
150 & 1000 & $6.42 \times 10^{-248}$ & $-247.19$ & ---      & 5 \\
\bottomrule
\end{tabular}
\end{table}

The observed positive-branch values
$\log_{10}|\lambda_{\min}^{\mathrm{even}}(N)|$ decrease monotonically
with $N$ at fixed $c=100$.  We note that the ordinary
Galerkin/Rayleigh--Ritz min-max monotonicity assumes a coercive
self-adjoint setup with an isolated ordered ground state; because the
negative-sign block at $c=100$ is a finite-$T$ artifact
(Section~\ref{sec:c100-doublets}), cutoff-free the smallest-positive
selection is the min-max ground state of the matrix; we nonetheless
report the monotonicity as an empirical observation, not as an
invocation of that theorem.  The successive first differences
$|\Delta_1|=56.28$, $|\Delta_2|=47.12$, $|\Delta_3|=39.37$ form a
geometric-looking sequence with consecutive ratios
$|\Delta_2/\Delta_1|=0.837$ and $|\Delta_3/\Delta_2|=0.836$, matching
to two decimal places --- evidence for a local geometric model of
the convergence sequence, which the three-point sweep
$\{100,150,200\}$ alone could not distinguish from a forced fit
(Section~\ref{sec:c100-aitken}).

\subsection{\texorpdfstring{Aitken $\Delta^2$ extrapolation}{Aitken Delta-squared extrapolation} and match to Connes 2026 \S6.4}
\label{sec:c100-aitken}

Define $x_N := \log_{10}|\lambda_N^{\mathrm{even}}(c=100)|$ on the
positive branch, with measured values
$x_{100}=-190.92$, $x_{150}=-247.19$, $x_{200}=-294.31$,
$x_{250}=-333.68$ (Table~\ref{tab:c100-N-sweep}).  First differences
$\Delta_i := x_{N_{i+1}} - x_{N_i}$ on the step grid $N\in\{100,150,200,250\}$
are $\Delta_1=-56.28$, $\Delta_2=-47.12$, $\Delta_3=-39.37$, with
consecutive ratios
\begin{equation}\label{eq:aitken-ratios}
  |\Delta_2/\Delta_1| = 0.837,\qquad
  |\Delta_3/\Delta_2| = 0.836,
\end{equation}
matching to two decimal places.  This consistency is the data
signature of a locally geometric convergence
$x_N = x_\infty + C\,r^{N}$; with four points and one shared ratio,
the geometric hypothesis is over-determined and therefore testable
on the data (in contrast to the three-point case, where any single
ratio fits both gaps by construction).

The Aitken-$\Delta^2$ acceleration \cite{Aitken1926} applied to two
consecutive overlapping triples yields two independent extrapolations:
\begin{equation}\label{eq:aitken-c100}
\begin{aligned}
\text{from } \{x_{100},x_{150},x_{200}\}:\quad
  x_\infty &\;\approx\; x_{100} - \frac{(\Delta_1)^2}{\Delta_2-\Delta_1}
            \;=\; -190.92 - \frac{(-56.28)^2}{9.16}
            \;\approx\; -536.76, \\[2pt]
\text{from } \{x_{150},x_{200},x_{250}\}:\quad
  x_\infty &\;\approx\; x_{150} - \frac{(\Delta_2)^2}{\Delta_3-\Delta_2}
            \;=\; -247.19 - \frac{(-47.12)^2}{7.75}
            \;\approx\; -533.70.
\end{aligned}
\end{equation}

The Aitken anchors in equation~\eqref{eq:aitken-c100} are computed
from the full-precision $\log_{10}|\lambda_{\min}^{\rm even}|$ values
recorded in the released JSON; the rounded table values in
Table~\ref{tab:c100-N-sweep} reproduce them only to the displayed
precision.

\begin{figure}[!htbp]
\centering
\includegraphics[width=\textwidth]{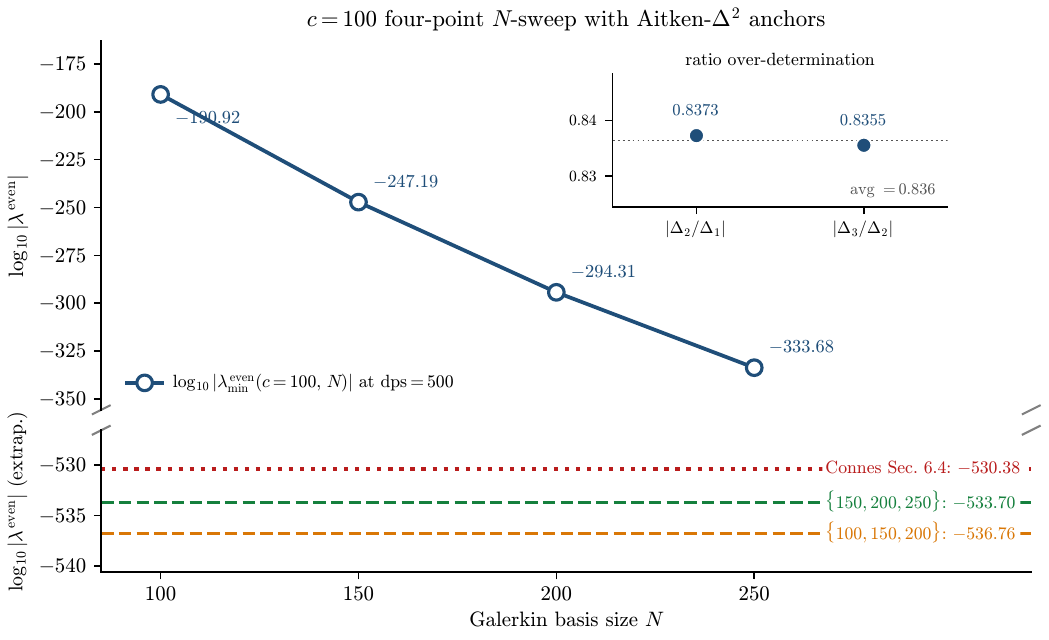}
\caption{Broken-axis $N$-sweep at $c=100$, $\dps=500$
  (smallest-positive even-sector eigenvalue).  \emph{Top panel:} the
  four measured data points $\log_{10}|\lambda_{\min}^{\rm even}|
  = -190.92, -247.19, -294.31, -333.68$ at $N=100, 150, 200, 250$,
  shown at their full ${\sim}140$-OOM range.  \emph{Bottom panel}
  (note the y-axis break and the change of scale): the two
  consecutive Aitken-$\Delta^2$ extrapolations
  ($-536.76$, $-533.70$) and the Connes 2026 \S6.4 heuristic continuum
  prediction ($-530.38$), at their natural ${\sim}10$-OOM scale.
  The trend with $N$ is monotone toward the Connes prediction; the
  deeper-anchored Aitken triple $(N\in\{150,200,250\})$ is within
  $3.32$~OOM of the heuristic.  \emph{Inset (top panel):} consecutive
  first-difference ratios $|\Delta_{i+1}/\Delta_i|$ from
  equation~\eqref{eq:aitken-ratios}.  Both ratios lie within $0.002$
  of the dotted reference line at $0.836$ (the average), evidence for
  a local geometric model of the convergence sequence (not a proof
  that alternatives are ruled out; see
  Section~\ref{sec:c100-aitken} discussion).}
\label{fig:c100-aitken}
\end{figure}

The two Aitken anchors agree to $|{-536.76}-({-533.70})| \approx
3.1$ orders of magnitude; the second (using the deeper triple) is
closer to the Connes 2026 \S6.4 heuristic prediction of
$-530.38$ in equation~\eqref{eq:connes-100-prediction}.  The gap
of the first anchor to the Connes prediction is $6.39$~OOM and the
gap of the second anchor is $3.32$~OOM, out of a magnitude range
$|x_\infty|\sim 530$.  The trend with $N$ is monotone toward the
prediction (Aitken anchored at $N=100$: $-536.76$; anchored at
$N=150$: $-533.70$; Connes 2026: $-530.38$).  Figure~\ref{fig:c100-aitken}
summarizes the four-point sweep, the two consecutive Aitken anchors,
and the Connes 2026 \S6.4 reference together with the
consecutive-ratio consistency check.

We make no claim that this constitutes a proof of the asymptotic.
Connes 2026 \S6.4 itself is a heuristic statement based on the
empirical agreement (their Figure~1) between $1-\chi_2(\lambda)$
and $\varepsilon(\lambda)$ for $\lambda\leq14$; the prefactor
$(2^{14}/3)\,\sqrt{2}\,\pi^{5}$ is derived from Fuchs's prolate-spheroidal
eigenvalue theorem and is therefore rigorous, while the
\emph{identification} of $\varepsilon$ with $1-\chi_2$ at large
$\lambda$ is the heuristic step.  On our side, the consistency of
the two ratios in equation~\eqref{eq:aitken-ratios} is evidence for
a local geometric model of the convergence sequence, but four points
do not rule out alternative convergence-model fits (stretched-exponential
$a_N = a_\infty + Ce^{-aN^\beta}$ with $\beta\neq 1$, or polynomial
$1/N^p$); such alternatives fit the same four points less tightly than
the geometric form, but are not eliminated, and they yield
extrapolated limits outside the Aitken band.  The extrapolation
remains an inference about the $N\to\infty$ limit of an empirically
distinguished positive branch (Remark~\ref{rem:ccm-hyp-c100}).  We
report the consistency with the Connes 2026 \S6.4 heuristic; we do
not promote it to a statistical test.

\subsection{\texorpdfstring{$\gamma_k$}{gamma\_k} extraction up to 307--329 matching digits}
\label{sec:c100-gamma}

The eigenvector corresponding to $\lambda_{\min}^{\mathrm{even}}(c=100)$
on the empirically distinguished positive branch (see
Remark~\ref{rem:ccm-hyp-c100}) extracts $\gamma_k$ via the
$F_{\mathrm{even}}$ test function defined in \S\ref{sec:eigensolver},
equivalently (under the unitary equivalence of \S\ref{sec:spectral-triple})
the zero set of $\widehat{\xi}_N(z)$.

\begin{table}[ht]
\centering
\caption{Extraction of $\gamma_k$ at $c=100$ from the smallest-positive
eigenvector.  Matching-digit counts $\lfloor -\log_{10}|\gamma_k -
\mathtt{mpmath.zetazero}(k).\mathtt{imag}|\rfloor$, with findroot
tolerance $10^{-380}$ at the working precision shown.  Detected
$\gamma_k$ values are deposited in the ancillary file
listed in Section~\ref{sec:code} for independent verification.}
\label{tab:c100-gamma}
\begin{tabular}{@{}rrrr@{}}
\toprule
$k$ & $N{=}150$, $\dps{=}500$ & $N{=}150$, $\dps{=}1000$ & $N{=}250$, $\dps{=}500$ \\
\midrule
1  & 128 & 242 & \textbf{329} \\
2  & 122 & 239 & 325 \\
3  & 126 & 236 & 323 \\
4  & 117 & 233 & 320 \\
5  & 128 & 231 & 318 \\
6  & 119 & 228 & 316 \\
7  & 118 & 226 & 313 \\
8  & 128 & 224 & 312 \\
9  & 111 & 221 & 309 \\
10 & 126 & 219 & 307 \\
\bottomrule
\end{tabular}
\end{table}

At $c=100$, $N=150$, $\dps=1000$ the first ten Riemann zeros are
reproduced to $219$--$242$ matching digits; at $N=250$, $\dps=500$
the same extraction reaches $307$--$329$ matching digits.  For
reference, \CCM\ \cite{CCM2025} Section~6 reports $\gamma_1$ matching
to $\sim55$~digits at $c=13$, $N=120$; our $c=67$, $N=100$, $\dps=200$
datum, whose error is $1.478\times10^{-168}$, gives $167$ matching
digits under the definition above (see
Table~\ref{tab:appendix-results}).

\begin{figure}[!htbp]
\centering
\includegraphics[width=0.85\textwidth]{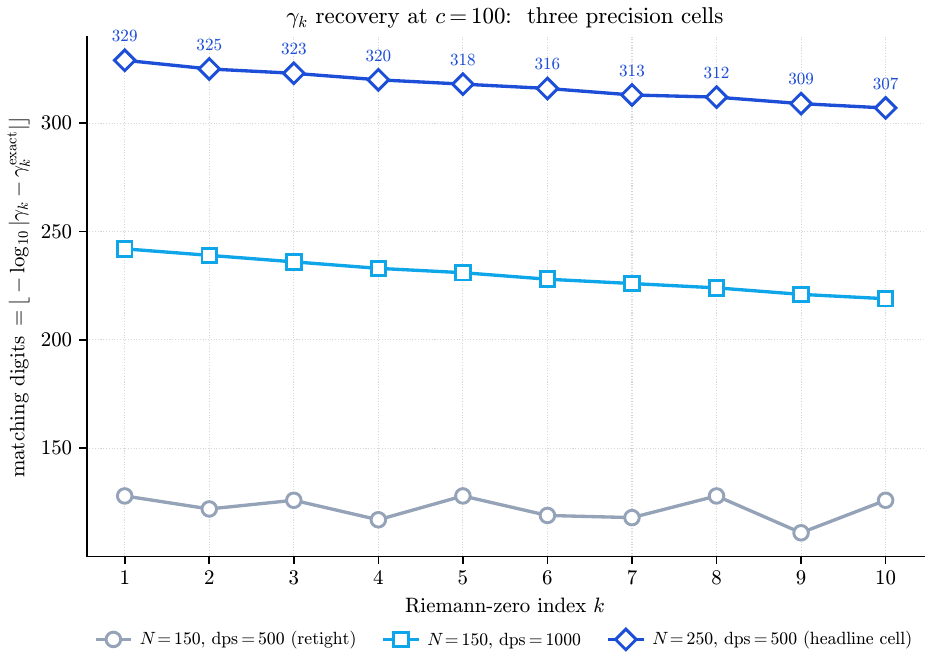}
\caption{Matching-digit recovery of $\gamma_1,\ldots,\gamma_{10}$
  at $c=100$ from three precision cells.  $N=150$ at $\dps=500$
  (retight-tolerance baseline) gives $111$--$128$ digits;
  $N=150$ at $\dps=1000$ gives $219$--$242$ digits (precision-doubling
  at fixed $N$); $N=250$ at $\dps=500$ gives $307$--$329$ digits
  (the headline cell of Table~\ref{tab:c100-gamma}).  Each annotation
  on the headline curve is the floor of $-\log_{10}|\gamma_k -
  \mathrm{zetazero}(k)|$ recomputed independently from the deposited
  $400$-digit reference strings.}
\label{fig:c100-gamma-digits}
\end{figure}

The doubling in $\dps$ ($500\to 1000$) at fixed $N=150$ adds
$93$--$117$ matching digits across $\gamma_1$ through $\gamma_{10}$;
the increment in $N$ ($150\to 250$) at $\dps=500$ adds $181$--$203$
matching digits at the same $\gamma_k$.  The $N$-direction therefore
gives more matching digits per compute-time unit than the $\dps$
direction at $c=100$: increasing $N$ both deepens the Galerkin
truncation \emph{and} provides cheaper accuracy than the $\dps$
quadrature cost.  This identifies the Galerkin truncation $N$ as
the binding accuracy parameter at $c=100$, with $\dps$ secondary
once it is comfortably above the backward-error floor.
Figure~\ref{fig:c100-gamma-digits} visualizes the three precision
cells; the gap between the $N=150,\dps=500$ and $N=150,\dps=1000$
curves isolates the precision-floor effect at fixed $N$, and the gap
between $N=150,\dps=1000$ and $N=250,\dps=500$ isolates the
$N$-extension effect at fixed working precision.

\subsection{The dps-stable negative block at \texorpdfstring{$c=100$}{c=100} is a finite-cutoff artifact}
\label{sec:c100-doublets}

At $c=100$, $N=150$, $T=800$, $\dps=500$, the even-sector matrix
$Q_c^{\mathrm{even}}$ has five eigenvalues with negative sign at
$\log_{10}|e|\in\{-40.66, -59.75, -96.75, -165.15, -223.69\}$.  The
$\dps=1000$ retest reproduces all five at log-magnitudes identical
to two decimal places, indicating they are not precision-floor noise
(the topmost negative at $10^{-41}$ lies hundreds of orders of
magnitude above the $\dps=1000$ working-precision floor).

Continuum positivity of $QW_\lambda$ for all $\lambda > 1$ is
equivalent to the Riemann Hypothesis (Connes 2026~\cite{Connes2026}
\S4.1) and is proved only for small $\lambda$; we therefore do not
assume continuum positivity at $\lambda = \sqrt{100}$.  The block is
$\dps$-stable but not $T$-stable: it is an artifact of truncating the
archimedean integral at the finite cutoff $T=800$.  Raising the cutoff
does not simply shrink the block, it rearranges it: at this cell the
five negatives at $T=800$ are replaced at $T=1200$ by three at entirely
different magnitudes ($\log_{10}|e|\in\{-71.73,-101.66,-176.21\}$), so
the count declines without reaching zero at any cutoff we tested.  A
cutoff-free evaluation
of the archimedean entries, certified by a rigorous interval
$LDL^{\mathsf T}$ factorization in Arb, leaves the $c=100$ even sector
with no negative eigenvalues at $N=100$, $150$ and $200$.  The correct
diagnostic of a deep-spectrum value is therefore agreement between two
values of $T$, not stability in $\dps$.  At
$T=800$ the archimedean entries carry a truncation error far larger
than the reported negative magnitudes, so those signs are not
informative about the operator.

The count of negative-sign eigenvalues grows with $N$ at fixed
$c=100$ ($\{3, 5, 8, 11\}$ for $N\in\{100,150,200,250\}$, see
Table~\ref{tab:c100-N-sweep}), consistent with this mechanism: at
fixed $T$ the higher-index basis functions added as $N$ grows carry
archimedean entries that require a larger cutoff to resolve, so more
of the lowest eigenvalues fall below the finite-$T$ entry-error
floor.  We do not infer anything about the continuum operator from
these signs.  This is distinct from spectral pollution proper
(spurious eigenvalues
arising from gaps in the essential spectrum), studied for
self-adjoint operators in
\cite{DaviesPlum2004,LevitinShargorodsky2004,Parlett1998}; in the
modern operator-theoretic framework spectral inclusion holds for the
self-adjoint Galerkin sequence, but freedom from spectral pollution
has not been proved, and we do not claim it.

This finite-basis observation does not affect the smallest
\emph{positive} eigenvalue $\lambda_{\min}^{\mathrm{even}}=
6.42\times10^{-248}$ at $N=150$ reported in
Table~\ref{tab:c100-N-sweep}.  That value is the corroborated
positive branch on which the $\gamma_1$ extraction of
Section~\ref{sec:c100-gamma} operates and the quantity participating
in the Aitken extrapolation of Section~\ref{sec:c100-aitken}.  We
flag this finite-$T$ artifact here for completeness; the same
mechanism accounts for the apparent $\chi_3$ positivity breakdown observed
at $c=23,29$ in Section~\ref{sec:chi3-odd}, where the even-sector
negative is likewise removed once $T$ is increased, and in both cases
we report the positive-branch observable without claiming a
theorem-level interpretation of the negative block.

\subsection{Reframing of the \texorpdfstring{$N=100$}{N=100} empirical fit}
\label{sec:c100-reframe}

Section~\ref{sec:extension} established the empirical fit
$|\log_{10}\lambda_{\min}^{\mathrm{even}}(c)|\approx 13.24\,c^{0.634}$
on the $c\leq67$, $N=100$ data.  Extrapolating that fit to $c=100$
predicts $|\log_{10}\lambda_{\min}|\approx 245.4$.  Our $N=200$
observation at $c=100$ is $294.31$, exceeding the prediction by
$49$~OOM in the direction of \emph{faster} decay.  Refitting
$|\log_{10}\lambda|=A\,c^B$ on $\{c\leq67\}\cup\{c=100\}$ data degrades
the RMS residual from $9.0$~OOM (in-sample at $N=100$) to $13.9$~OOM
once the $c=100$ point is included --- a clean falsification of the
pure-power-law functional form.

The natural interpretation, consistent with the data, is that
$|\log_{10}\lambda|\sim 13.24\,c^{0.634}$ describes the
\emph{finite-$N=100$ rate} for $c\leq 67$, not the continuum rate.  The
continuum rate, per Connes 2026 \S6.4 and the Aitken extrapolation in
Section~\ref{sec:c100-aitken}, is super-exponential in $c$ ($\sim
e^{-4\pi c+9\log c/2}$).

This interpretation is independently corroborated by an additional
measurement at the deepest in-sample cutoff $c=67$ at the same
extended basis $N=150$.  Re-running the operator at $c=67$, $T=800$,
$N=150$, $\dps=500$ gives
$\lambda_{\min}^{\mathrm{even}}(c{=}67,\,N{=}150) =
5.329\times 10^{-219}$, i.e.\ $\log_{10}|\lambda_{\min}| = -218.27$.
This is $46$~OOM smaller than the $N=100$ value of
$\log_{10}|\lambda_{\min}(c{=}67, N{=}100)| = -172.10$ reported in
Section~\ref{sec:extension}.  The N-dependence at $c=67$
($-172.10\to -218.27$ over $N=100\to 150$, a drop of $46$~OOM) is
qualitatively the same as at $c=100$
($-190.92\to -247.19$, a drop of $56$~OOM).  The $N=100$ data of
Section~\ref{sec:extension} across the entire $c\leq67$ sweep are,
like the $N=100$ datum at $c=100$, Galerkin upper bounds rather than
near-continuum values.
This is the structural mechanism underlying the failure of the
$13.24\,c^{0.634}$ extrapolation at $c=100$: the empirical fit
captures a sequence of upper bounds rather than the limit they
approximate from above.

The $N=100$ finite-$N$ fit of Section~\ref{sec:extension} and
Connes' continuum asymptotic measure different quantities; the
present work makes that distinction explicit, with both quantities
now empirically tested at $c=100$.

The four observations of this section --- the two consecutive Aitken
anchors at $-536.76$ and $-533.70$ approaching the Connes 2026 \S6.4
heuristic of $-530.38$ with gaps $6.39$ and $3.32$ OOM respectively,
the $307$--$329$ digit reproduction of $\gamma_1$--$\gamma_{10}$ from
the $N=250$ smallest-positive eigenvector at $\dps=500$ (with the
$N=150, \dps=1000$ retest giving $219$--$242$ digits at the same
$\gamma_k$), the finite-$T$ negative-sign artifact block (counts
$\{3,5,8,11\}$ at $T=800$, removed at larger $T$), and the falsification
of the finite-$N$ power-law extrapolation by $49$~OOM at $N=200$ ---
are the substantive new content of this work relative to the $c\leq67$
data of Section~\ref{sec:extension}.  They constrain theoretical work on
the convergence question through structural quantities (Sobolev
exponent, eigenvector universality) which
the subsequent Section~\ref{sec:structural} examines in detail.

\section{Precision-Floor Analysis and the
  \texorpdfstring{$\dps=150$}{dps=150} Retest}
\label{sec:precision}

\subsection{The backward-error precision floor}

At $\dps=80$, the unit roundoff is $\varepsilon\approx10^{-80}$.  The
backward-error bound for symmetric eigenvalue computation gives
\begin{equation}\label{eq:backward-error}
  |\delta\lambda| \leq \varepsilon\cdot\|Q\|_2,
\end{equation}
where $\|Q\|_2$ is the spectral norm of the Galerkin matrix.  Direct
computation gives $\|Q\|_2\approx5.7$--$6.0$ across the cutoffs
tested.  The backward-error floor is therefore
\begin{equation}\label{eq:floor-80}
  \varepsilon\cdot\|Q\|_2 \approx 6\times10^{-80}.
\end{equation}

At $c=13$, $\lambda_{\min}^{\mathrm{even}}=2.865\times10^{-59}$ is
21~orders of magnitude above this floor --- comfortably resolved.  At
$c=17$, $\lambda_{\min}^{\mathrm{even}}=2.030\times10^{-80}$ is
\emph{at or below} the floor.  At $c=19$,
$\lambda_{\min}^{\mathrm{even}}=1.265\times10^{-90}$ is firmly below
the floor at $\dps=80$.

\subsection{Pre-registered five-bin classification rule}

Before the $\dps=150$ retest was executed, a five-bin classification
rule was frozen based on the quantity
\[
  \Delta_{17\text{-}19}^{150}
  = \log_{10}|\gamma_1\text{ error at }c=19,\;\dps=150|
  - \log_{10}|\gamma_1\text{ error at }c=17,\;\dps=150|.
\]

\begin{enumerate}
\item \textbf{5A confirmed} ($\Delta<-15$): The $c=19$ datum was
  entirely an artifact of the precision floor.
\item \textbf{Ambiguous partial floor} ($-15\leq\Delta<-5.29$): Some
  floor contamination, but the step pattern is not fully explained.
\item \textbf{5B live} ($-5.29\leq\Delta<-1$): The non-smooth step
  pattern is largely genuine.
\item \textbf{5B strong} ($-1\leq\Delta<0$): Near-stalling at $c=19$
  is real.
\item \textbf{Reversal} ($\Delta\geq0$): The error increases from
  $c=17$ to $c=19$.
\end{enumerate}

\subsection{Retest results}

The $\dps=150$ retest completed in approximately 118~minutes of
wall-clock time.

\begin{table}[ht]
\centering
\caption{$\dps=80$ versus $\dps=150$ comparison.}
\label{tab:dps-retest}
\begin{tabular}{@{}llll@{}}
\toprule
Quantity & $\dps=80$ & $\dps=150$ & Shift \\
\midrule
$\lambda_{\min}$ at $c=13$
  & $2.077\times10^{-59}$ & $2.077\times10^{-59}$
  & $<10^{-59}$ shift \\
$|\gamma_1\text{ error}|$ at $c=13$
  & $1.455\times10^{-55}$ & $1.455\times10^{-55}$
  & $<10^{-59}$ shift \\
$\lambda_{\min}$ at $c=17$
  & $1.215\times10^{-80}$ & $1.401\times10^{-80}$
  & $+15\%$ \\
$|\gamma_1\text{ error}|$ at $c=17$
  & $1.117\times10^{-76}$ & $1.128\times10^{-76}$
  & $+1\%$ \\
$\lambda_{\min}$ at $c=19$
  & $1.327\times10^{-81}$ & $8.777\times10^{-91}$
  & 9 OOM drop \\
$|\gamma_1\text{ error}|$ at $c=19$
  & $5.809\times10^{-80}$ & $7.426\times10^{-87}$
  & 7 OOM drop \\
\bottomrule
\end{tabular}
\end{table}

The $c=19$ cell was substantially contaminated at $\dps=80$: the
eigenvalue dropped by 9~orders of magnitude and
$|\gamma_1\text{ error}|$ dropped by 7~orders of magnitude under the
precision increase.  The $c=13$ cell is unaffected (agreement to all
reported digits).

\subsection{Classifier verdict}

At $\dps=150$:
\begin{equation}\label{eq:classifier}
  \Delta_{17\text{-}19}^{150}
  = \log_{10}(7.426\times10^{-87}) - \log_{10}(1.128\times10^{-76})
  = -10.18.
\end{equation}
This falls in \textbf{Bin~2: Ambiguous partial floor}
($-15\leq-10.18<-5.29$).  The classifier fires without post-hoc
adjustment.

\textbf{Interpretation.}  The $c=19$ $\dps=80$ datum was substantially
contaminated (the 7-OOM drop in $|\gamma_1\text{ error}|$ is direct
evidence), but the step from $c=17$ to $c=19$ at $\dps=150$
($\Delta=-10.18$) is still smaller than the step from $c=13$ to $c=17$
($\Delta=-21.11$).

\subsection{The \texorpdfstring{$|\gamma_1\text{ err}|/\lambda_{\min}$}%
  {|gamma1 err|/lambda\_min} coupling}

\begin{figure}[!htbp]
\centering
\includegraphics[width=0.85\textwidth]{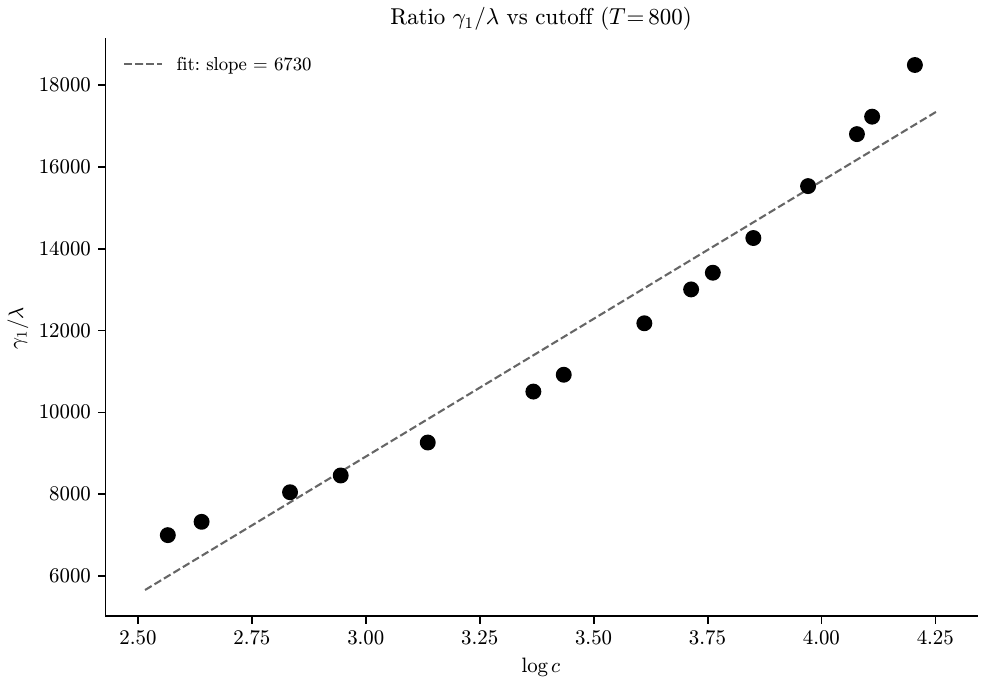}
\caption{The ratio $|\gamma_1\text{ err}|/\lambda_{\min}$ as a function of
  $\log c$, showing slow monotone growth consistent with standard
  spectral approximation theory.  The dashed line is the
  linear fit $C(c)\approx6730\cdot\log c - 11268$.}
\label{fig:ratio}
\end{figure}

The ratio $|\gamma_1\text{ err}|/\lambda_{\min}$ ranges from approximately 6999 at
$c=13$ to 18489 at $c=67$ (at $T=800$, $\dps=200$), exhibiting a
slow upward drift with increasing~$c$.  This ratio is an empirical
observation consistent with standard Babu\v{s}ka--Osborn /
Rayleigh--Ritz spectral approximation theory
\cite{BabuskaOsborn1991}; the prefactor $C(c)$
grows slowly, approximately linearly in $\log(c)$
(see Figure~\ref{fig:ratio}).

A linear fit to the 15-point $T=800$ data gives
\[
  C(c) = |\gamma_1\text{ err}|/\lambda_{\min}
  \approx 6730\cdot\log c - 11268,
  \qquad R^2 = 0.956.
\]

The coupling is also visible in cross-cutoff log-ratio scaling:

\begin{table}[ht]
\centering
\caption{Cross-cutoff log-ratio coupling ($T=800$, $\dps=200$).}
\label{tab:log-ratio}
\begin{tabular}{@{}lccc@{}}
\toprule
Pair & $\log_{10}(\lambda\text{ ratio})$
     & $\log_{10}(\gamma_1\text{ ratio})$
     & Ratio of logs \\
\midrule
$(c=13,\;c=17)$ & 21.150 & 21.089 & 1.003 \\
$(c=17,\;c=19)$ & 10.206 & 10.184 & 1.002 \\
$(c=13,\;c=19)$ & 31.355 & 31.273 & 1.003 \\
$(c=13,\;c=43)$ & 89.943 & 89.725 & 1.002 \\
$(c=13,\;c=67)$ & 113.555 & 113.131 & 1.004 \\
\bottomrule
\end{tabular}
\end{table}

The eigenvalue and first-zero error scale in near-perfect lockstep in
log space: the ratio of logs is within 0.4\% of unity across the full
113-OOM span.  This log-space coupling
$\log(\gamma_1\text{ ratio})/\log(\lambda\text{ ratio})\approx1.003$
held across all measured pairs.

For reference, the $c=19$ $\dps=80$ ratio was 44, an obvious artifact;
at $\dps=150$ it restores to 8461, consistent with floor contamination.

\subsection{Eigenvector robustness to eigenvalue contamination}

An initially surprising feature of the retest is that
$|\gamma_1\text{ error}|$ at $c=17$ shifted by only 1\% under the
$\dps=80\to150$ transition, even though $\lambda_{\min}$ shifted by
15\%.  The explanation follows from standard eigenvector perturbation
theory (Davis--Kahan): the eigenvector perturbation depends on the
spectral gap, which
at $c=17$ is approximately $10^{8}$.  The eigenvector is therefore
well-conditioned even when the eigenvalue is at the precision floor, and
$\gamma_1$ --- computed from the eigenvector's Fourier transform ---
inherits this stability.

This observation has a practical consequence:
$|\gamma_1\text{ error}|$ is a more reliable observable than
$\lambda_{\min}$ when operating near the precision floor.

\section{Structural Analysis}\label{sec:structural}

\subsection{Empirical Galerkin convergence exponent (Sobolev interpretation)}
\label{sec:sobolev-scaling-emp}

Standard Galerkin approximation theory
\cite{BabuskaOsborn1991}
predicts $|\lambda_{\mathrm{exact}}-\lambda_N|\sim C\cdot N^{-2s}$
where~$s$ is the Sobolev regularity exponent of the target
eigenfunction.  The $N$-convergence data at $c=23$
(Table~\ref{tab:N-conv}) allow a direct \emph{empirical} measurement
of the convergence-rate exponent, which under the Babu\v{s}ka--Osborn
identification \emph{would} equal $s$ if the standard hypotheses
hold (smooth coercive bilinear form, isolated simple eigenvalue,
target eigenfunction in $H^s$).  We do not establish those hypotheses
analytically here; the values $s(c)$ below should be read as
\emph{empirical Galerkin-convergence rates} on the measured $N$-range,
not a proof that $\eta_c \in H^s$.

Linear regression of $\log_{10}|\lambda_{\min}^{\mathrm{even}}|$
against $\log_{10}N$ for $N\in\{40,60,80\}$ (the pre-saturation
regime) yields slope $a=-92.28$ with $R^2=0.99997$.  Since
$|\lambda|\sim N^a$ and $a=-2s$, this gives
\[
  s \approx 46.1 \quad (\text{from }\lambda_{\min}), \qquad
  s \approx 46.5 \quad (\text{from }|\gamma_1\text{ error}|),
\]
both with $R^2>0.9999$.  The two independent estimates agree to
within~0.7\%.

\begin{figure}[!htbp]
\centering
\includegraphics[width=0.85\textwidth]{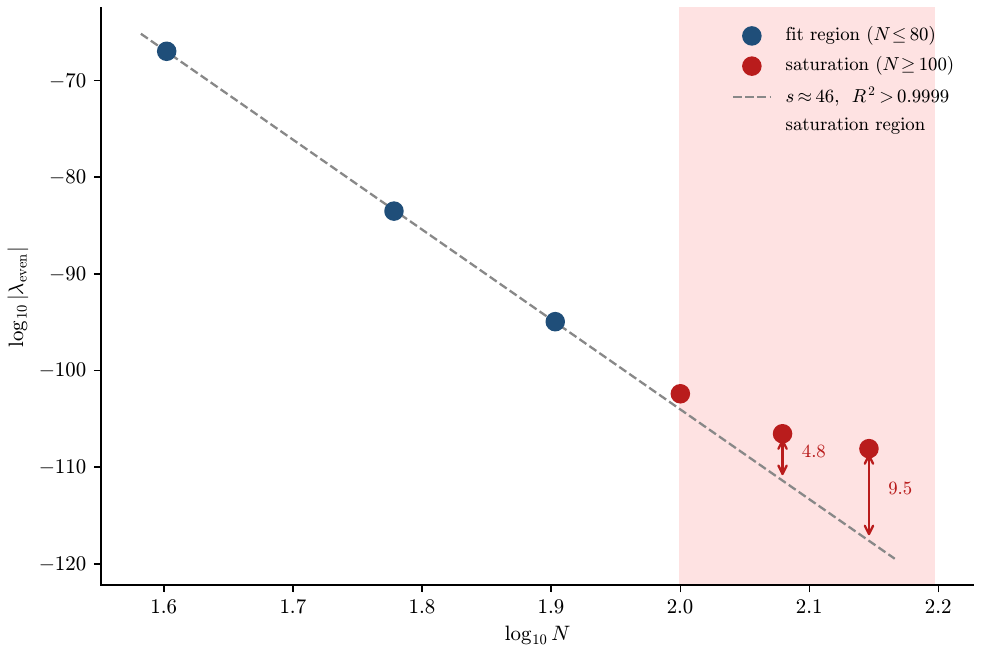}
\caption{Sobolev regularity measurement at $c=23$.  The power-law
  fit through $N=40,60,80$ yields $s\approx46$ with $R^2>0.9999$.
  Saturation at $N\geq100$ (shaded) indicates the $\dps=150$ precision
  floor has been reached.}
\label{fig:sobolev}
\end{figure}

At $N\geq100$, the data fall above the power-law prediction by
1.6--9.4~decades (Figure~\ref{fig:sobolev}), indicating that the
$\dps=150$ precision floor has been reached and the Galerkin
truncation error is no longer the dominant error source.  The
crossover occurs at $N\approx100$, consistent with the observation
that $k_{\mathrm{eff}}\approx46$ modes carry the ground-state
eigenvector.

Under the Babu\v{s}ka--Osborn identification, this would correspond
to $H^{46}([0,L])$-type regularity for the CvS ground-state
eigenvector at $c=23$; in this paper we use $s(c)$ only as an
empirical Galerkin-convergence exponent on the measured $N$-range,
not as an established Sobolev statement.  Repeating the
$N$-convergence measurement at five additional cutoffs reveals that
the empirical exponent is \emph{not} a universal constant --- it
grows with~$c$:

\begin{table}[ht]
\centering
\caption{Sobolev regularity exponent $s(c)$ at six cutoffs,
  each measured from the pre-saturation regime
  $N\in\{40,60,80\}$ at $T=400$, $\dps=150$.  The archimedean cutoff is
  $T=400$, inherited from the per-cutoff $N$-convergence ladders
  (Table~\ref{tab:N-conv}), not the $T=800$ of the fifteen-cutoff sweep.}
\label{tab:sobolev-scaling}
\begin{tabular}{@{}rrr@{}}
\toprule
$c$ & $s(c)$ & $R^2$ \\
\midrule
13 & 9.4  & 0.87 \\
17 & 27.8 & 0.99 \\
23 & 46.1 & 1.00 \\
29 & 57.6 & 1.00 \\
37 & 68.8 & 1.00 \\
43 & 75.0 & 1.00 \\
\bottomrule
\end{tabular}
\end{table}

A linear fit in $\log c$ gives
\[
  s(c) \approx 55\cdot\log c - 128,
  \qquad R^2 = 0.992.
\]
The empirical exponent $s(c)$ increases from $\approx 9$ at $c=13$ to
$\approx 75$ at $c=43$ --- under the Babu\v{s}ka--Osborn
identification, this \emph{would} correspond to a substantial
increase in Sobolev regularity, but again, we use $s(c)$ here as a
Galerkin-rate measurement, not as a regularity proof.  The scaling
has a direct implication for the observed convergence rate: at
$c=43$ with $s(c)\approx75$, the rate $N^{-2s(c)} = N^{-150}$ is so
rapid that $N=100$ basis functions reach the precision floor.

The low $R^2=0.87$ at $c=13$ reflects the fact that even $N=60$ and
$N=80$ are approaching saturation at $\dps=150$ for this cutoff; the
pre-saturation regime is narrower when $s$ is small.  For $c\geq17$,
all fits have $R^2>0.98$.

\subsection{Sobolev--Paley--Wiener connection: a proposed mechanism, and
its refutation}
\label{sec:sobolev-pw}

The Sobolev scaling $s(c)\approx 55\cdot\log c - 128$ appeared to admit a
mechanism via the classical Paley--Wiener theorem.  The proposal is stated
here as it was, followed by the measurement that refutes it.  The CvS
kernel involves the digamma function $\Psi(1/4 + i\tau/2)$, which is
analytic in a strip of width $\sigma=1/2$ around the real axis.  For a
Galerkin method on $[0,T]$, the Paley--Wiener theorem predicts
\[
  s = \frac{\sigma_{\mathrm{eff}}\cdot T}{2\pi},
\]
where $\sigma_{\mathrm{eff}}$ is the effective analyticity strip width
of the full kernel.  The bare digamma strip width ($\sigma=1/2$) is
independent of~$c$, but the prime-sum kernel --- a sum over primes
$p\leq c$ weighted by $(\log p)^2/p$ --- smooths the kernel beyond
the analyticity strip of any single term.  By the explicit formula for
the prime-counting remainder, the fluctuations have amplitude
$\sim\sqrt{c}/\log c$, and the effective strip width grows as
$\sigma_{\mathrm{eff}}\sim A\cdot\log c$ for some constant~$A$.

Taking $T=800$ and the measured slope $55$, this gave
$A = 55\cdot 2\pi/800 = 0.432$, implying $\sigma_{\mathrm{eff}}\approx
0.43\cdot\log c$ at any given cutoff, with the testable prediction
$s(c,T) \sim (A/T_0)\cdot T\cdot\log c$.

\medskip
\noindent\textbf{The mechanism is withdrawn.}  Three things are wrong with the
paragraph above.  First, it misidentifies the interval.  The Galerkin space
here is spanned by the periodic modes on $[0,L]$ with $L=\log c$
(Section~\ref{sec:setup}, Section~\ref{sec:eigensolver}), and the
Babu\v{s}ka--Osborn reading invoked in Section~\ref{sec:sobolev-scaling-emp}
is posed on that same interval; the cutoff $T$ bounds the archimedean
$\tau$-integral of the quadrature (Section~\ref{sec:quadrature}) and plays no
part in the mode spacing $2\pi/L$.  Whatever exponent governs the Fourier
decay on $[0,L]$, it is not $\sigma_{\mathrm{eff}}T/2\pi$.  Second, the
arithmetic substitutes $T=800$ for a slope measured at $T=400$: the exponents
in Table~\ref{tab:sobolev-scaling} come from the per-cutoff $N$-convergence
ladders at $T=400$ (Table~\ref{tab:N-conv}), not from the fifteen-cutoff
sweep convention, so $A = 55\cdot 2\pi/800$ is not a well-formed reading of
the data.  Third, the relation as printed predicts $s\propto T$ at fixed~$c$,
and that prediction has now been tested directly.  Repeating the $c=23$
measurement on the same $N\in\{40,60,80\}$ grid at $\dps=150$ gives

\begin{center}
\begin{tabular}{@{}rl@{}}
\toprule
$T$ & $s$ \\
\midrule
$400$  & $46.140$ \\
$800$  & $46.031$ \\
$1600$ & $45.934$ \\
\bottomrule
\end{tabular}
\end{center}

\noindent A fourfold increase in~$T$ therefore moves the exponent by $0.206$,
where proportionality requires $s(1600)\approx 4\,s(400)\approx 184.6$.  A
single cutoff settles this, because the relation asserts proportionality at
fixed~$c$.  By the criterion pre-registered in Section~\ref{sec:future}, the
relation as printed is refuted, and the constant $A$ has no meaning as
derived.  What
survives is the measurement itself: the empirical scaling
$s(c)\approx 55\log c - 128$ is a fit to the measured exponents rather than a
consequence of the mechanism, and both it and the exponents are unaffected.
We have no mechanism for the $c$-dependence at present and record that as an
open question rather than replacing one plausible story with another.

The three-cutoff test was first carried out by M.~Osman, who reported it in
issue~2 of the repository listed in Section~\ref{sec:code}; the values above
are our independent recomputation, which reproduces his to the precision he
reported.

\subsection{Approximate eigenvector universality across cutoffs}

\begin{figure}[!htbp]
\centering
\includegraphics[width=0.85\textwidth]{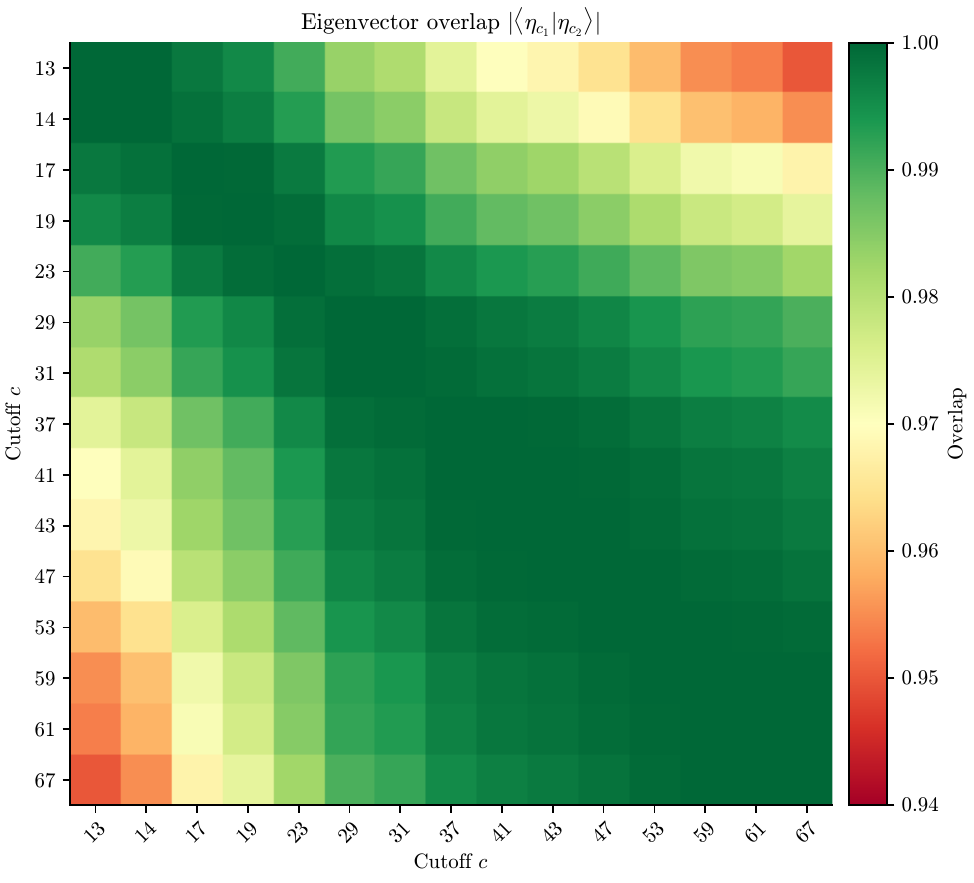}
\caption{Eigenvector overlap matrix $|\langle\eta_{c_1}|\eta_{c_2}\rangle|$
  across all fifteen cutoffs.  All 105~pairwise overlaps are at least
  $0.9498$ (the minimum, $0.94985$, occurs at the maximally-separated
  pair $(c=13,c=67)$; 104 of 105 pairs strictly exceed $0.950$),
  indicating approximate eigenvector universality despite eigenvalues
  differing by up to 113~orders of magnitude.}
\label{fig:overlap}
\end{figure}

The full $15\times15$ overlap matrix
$|\langle\eta_{c_1}|\eta_{c_2}\rangle|$ was computed from the saved
eigenvectors at all fifteen cutoffs.  All~105 pairwise overlaps are
at least $0.9498$ (Figure~\ref{fig:overlap}), with $104$ of $105$
strictly exceeding $0.950$ and the minimum
$0.94985$ attained at the most-separated pair $(c=13,c=67)$.
The closest-neighbor pairs $(c=41,43)$ and $(c=29,31)$ achieve overlaps
of $0.99997$ and $0.99992$ respectively.

\begin{figure}[!htbp]
\centering
\includegraphics[width=0.85\textwidth]{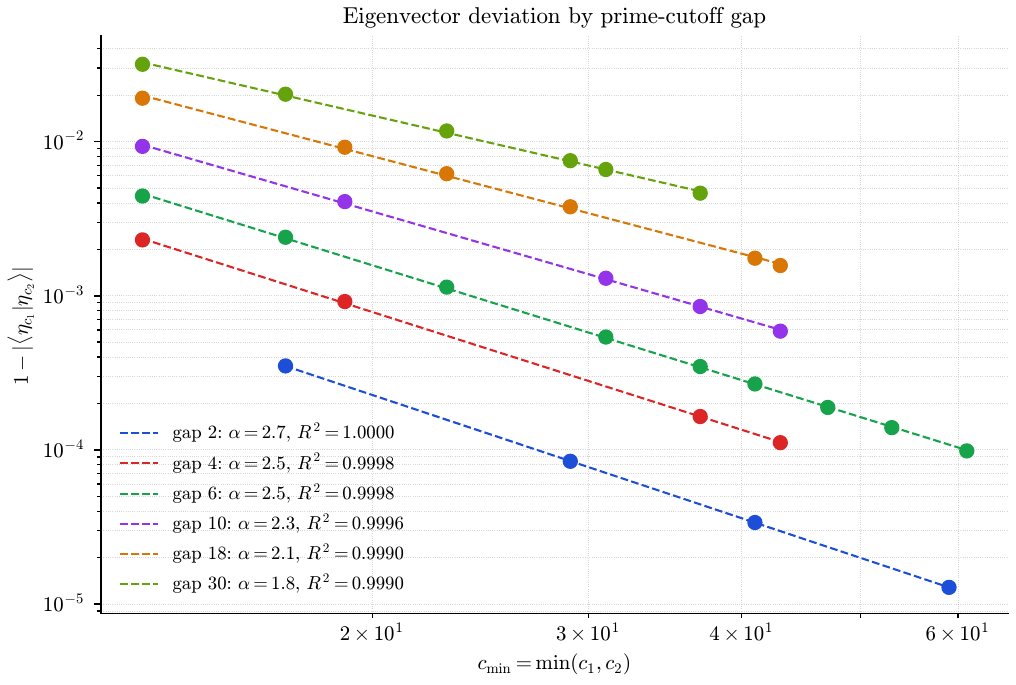}
\caption{Eigenvector deviation $1-|\langle\eta_{c_1}|\eta_{c_2}\rangle|$
  versus $c_{\min}=\min(c_1,c_2)$, grouped by prime-cutoff gap.
  At fixed gap, each series follows a clean power law
  $\sim c_{\min}^{-\alpha}$ with $R^2>0.999$.  The exponent ranges
  from $\alpha\approx2.7$ (gap~2) to $\alpha\approx1.8$ (gap~30).}
\label{fig:eigvec-universality}
\end{figure}

For fixed prime gap, the convergence rate follows a clean power law
(Figure~\ref{fig:eigvec-universality}):
\[
  1 - |\langle\eta_{c_1}|\eta_{c_2}\rangle|
  \sim A\cdot c_{\min}^{-\alpha}, \qquad c_{\min}=\min(c_1,c_2),
  \qquad \alpha \in [1.8,\,2.7], \quad R^2 > 0.999.
\]
This is consistent with standard rank-one perturbation theory: each
prime~$p$ contributes a rank-one correction of magnitude
$\Lambda(p)/\sqrt{p}\sim\log p/\sqrt{p}$ to the operator, and the
cumulative effect on the eigenvector scales as a power of~$c$.

The eigenvector profile broadens slowly with~$c$: at $c=13$,
components decay below $10^{-5}$ by $k=6$, while at $c=53$ this
threshold is not reached until $k\approx12$.  Despite this
broadening, the dominant components ($k=1$ through~5) are stable
across all cutoffs, and the effective support $\keff\approx46$ is
essentially $c$-invariant for $c\geq23$.

\subsection{Multi-zero convergence universality}
\label{sec:multi-zero}

The pickle files store the first ten detected zeros ($\gamma_1$
through~$\gamma_{10}$) at each cutoff.  This allows a test of whether
higher zeros converge at the same rate as~$\gamma_1$.

Linear regression of $\log_{10}|\gamma_k\text{ error}|$ against
$\log_{10}|\gamma_1\text{ error}|$ across cutoffs yields near-unity
slopes for all~$k$:

\begin{table}[ht]
\centering
\caption{Convergence-rate ratio $\text{slope}_k/\text{slope}_1$
  for the first five zeros.}
\label{tab:multi-zero}
\begin{tabular}{@{}lcl@{}}
\toprule
Zero & Rate ratio & Comment \\
\midrule
$\gamma_1$ & 1.000 & (reference) \\
$\gamma_2$ & 0.996 & 99.6\% of $\gamma_1$ rate \\
$\gamma_3$ & 0.993 & \\
$\gamma_4$ & 0.988 & \\
$\gamma_5$ & 0.985 & \\
$\gamma_{10}$ & 0.962 & 96.2\% of $\gamma_1$ rate \\
\bottomrule
\end{tabular}
\end{table}

\begin{figure}[!htbp]
\centering
\includegraphics[width=0.85\textwidth]{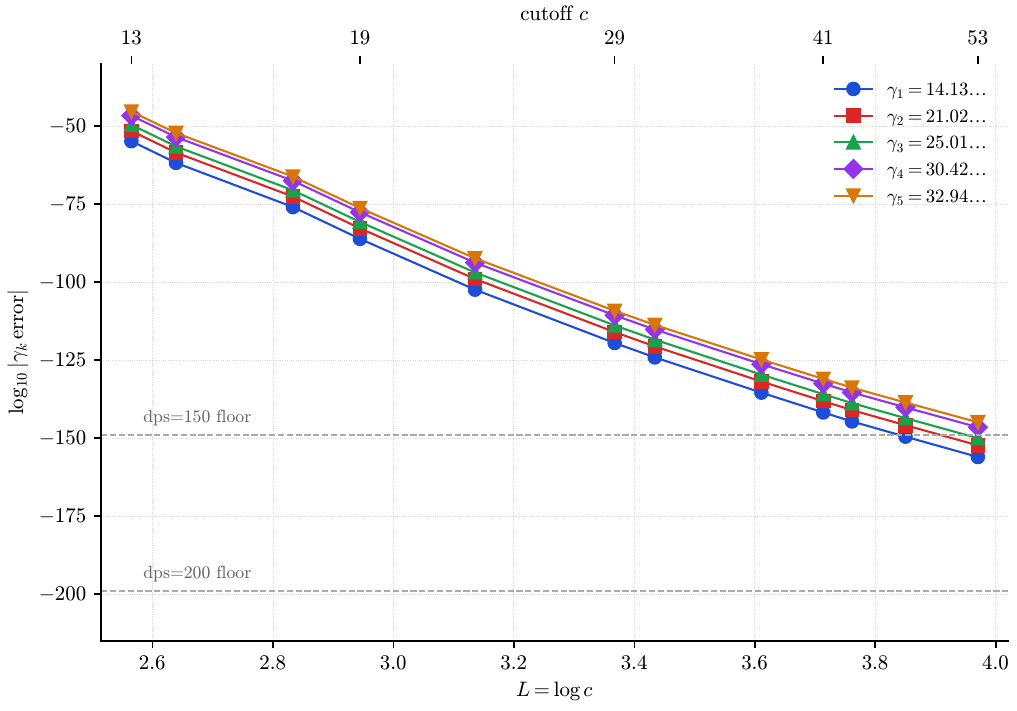}
\caption{Multi-zero convergence curves for $\gamma_1$ through
  $\gamma_5$ across all fifteen cutoffs (the curves for $\gamma_6$
  through $\gamma_{10}$ are visually indistinguishable from these and
  are omitted for legibility; per-zero rate ratios for $\gamma_1$
  through $\gamma_5$ and for $\gamma_{10}$ are tabulated in
  Table~\ref{tab:multi-zero}).  All ten zeros converge at rates
  within 3.8\% of each other, consistent with
  Heuristic~\ref{prop:multi-zero-universality}.}
\label{fig:multi-zero}
\end{figure}

All ten zeros converge at rates within 3.8\% of each other
(Figure~\ref{fig:multi-zero}).  The
$\gamma_2$ regression against~$\gamma_1$ has $R^2=0.9999998$.  This
near-universality means that the convergence question --- does
$|\gamma_1\text{ error}|\to0$ as $c\to\infty$? --- effectively
subsumes the corresponding question for all higher zeros.  If
$\gamma_1$ converges, the data strongly suggest that $\gamma_2$
through~$\gamma_{10}$ converge as well, at essentially the same rate.

The rate ratio decreases approximately linearly: each successive zero
converges about 0.42\% more slowly than~$\gamma_1$.  For $k\geq8$
at the largest cutoffs ($c\geq47$), the errors show anomalous
behavior consistent with finite-basis ($N=100$) truncation artifacts.

\begin{heuristic}[Empirical mechanism for multi-zero universality]
\label{prop:multi-zero-universality}
\emph{Working hypothesis:} let $Q_c$ be a family of self-adjoint
operators indexed by a cutoff~$c$, with ground-state eigenvectors
$v_c$ satisfying $\|v_c - v_{c'}\|_2 < \epsilon$ for all consecutive
cutoff pairs (a property we measure here; we do not assert
existence of a limit~$v_\infty$).  Let $\gamma_j[v]$ denote the
$j$-th zero of the Fourier--Mellin transform
$\mathrm{FM}(\gamma;v) = \sum_k v_k\,\phi_k(\gamma)$,
defined implicitly by $\mathrm{FM}(\gamma_j;v)=0$.  Suppose
$|\partial_\gamma\mathrm{FM}(\gamma_j;v)| \geq \delta > 0$ at the
tested zeros and cutoffs (verified numerically here; not established
analytically).  \emph{Under bound saturation}, that is, when the
perturbation $\delta v_c$ has a leading-order component in the
sensitivity direction at each $\gamma_j$, one expects the
relative-rate identity
\[
  \frac{|\gamma_j(c) - \gamma_j(\infty)|}{|\gamma_1(c) - \gamma_1(\infty)|}
  \;\lesssim\; 1 \;+\; O\!\left(\frac{\|\phi(\gamma_j)\|_2 / |\mathrm{FM}'(\gamma_j)|}
    {\|\phi(\gamma_1)\|_2 / |\mathrm{FM}'(\gamma_1)|} - 1\right)
\]
to hold to leading order in~$\epsilon$.
\end{heuristic}

\begin{proof}[Sketch (heuristic, not a formal proof)]
The Fourier--Mellin transform is linear in~$v$, hence Lipschitz with
constant $\|\phi(\gamma)\|_2$.  By the implicit function theorem,
\[
  \frac{d\gamma_j}{dv_k} = -\frac{\phi_k(\gamma_j)}{\partial_\gamma
  \mathrm{FM}(\gamma_j;v)}.
\]
The error $|\gamma_j(c) - \gamma_j(\infty)|$ is therefore bounded by
$\|\delta v_c\| \cdot \|\phi(\gamma_j)\|_2 / |\mathrm{FM}'(\gamma_j)|$.
Since $\delta v_c$ is the \emph{same vector} for all zeros~$j$, the
ratio of upper bounds reduces to the ratio of the sensitivity
coefficients $\|\phi(\gamma_j)\|_2/|\mathrm{FM}'(\gamma_j)|$.  The
measured Fourier--Mellin slopes vary by only 3.3\% across cutoffs,
so these coefficients are comparable across zeros, which would yield
the claimed universality \emph{provided} the bound saturates (i.e.,
the actual rate equals the upper bound to leading order).  This
saturation step is the heuristic content of the argument: it is
plausible at small $\epsilon$ but is not proved here.  Monotonicity
of $F_c$ in a neighborhood of each $\gamma_k$ is verified
numerically: $F'_c(\gamma_k)$ is strictly positive at all tested
cutoffs and zeros.
\end{proof}

\subsection{Eigenvector PCA: a two-dimensional manifold}
\label{sec:pca}

Principal component analysis of the fifteen normalized
eigenvectors --- stacked as rows of a $15\times101$ matrix (15~eigenvectors, each with 101~components) ---
reveals that they lie on an essentially two-dimensional manifold:

\begin{table}[ht]
\centering
\caption{PCA of the fifteen CvS ground-state eigenvectors.}
\label{tab:pca}
\begin{tabular}{@{}lrr@{}}
\toprule
Component & Variance explained & Cumulative \\
\midrule
PC1 & 98.73\% & 98.73\% \\
PC2 &  1.25\% & 99.98\% \\
PC3 &  0.02\% & 100.00\% \\
\bottomrule
\end{tabular}
\end{table}

Two principal components capture 99.98\% of the variance.  The
PC1 score is nearly perfectly correlated with $\log c$
($r=0.997$), meaning the eigenvector evolution from $c=13$ to
$c=67$ follows a one-parameter curve in the 101-dimensional
basis space.  Despite eigenvalues differing by 113~orders of
magnitude, the eigenvector's trajectory through Hilbert space
is effectively one-dimensional.

This is consistent with the existence of a limiting
eigenvector $\eta_\infty$: the eigenvectors are converging
along a fixed curve, not wandering through a high-dimensional
space.

\subsection{Frobenius-norm stability of the bulk spectrum}

The full CvS Galerkin matrix at $N=100$ has dimension~101.  While
the even-sector $\lambda_{\min}$ varies by 113~orders of magnitude
across cutoffs,
the Frobenius norm $\|Q_{\mathrm{even}}\|_F =
(\sum_i\lambda_i^2)^{1/2}$ is stable: it ranges from
30.74 ($c=31$) to 31.44 ($c=41$), a variation of only~2.3\%.
This stability reflects the dominance of the bulk eigenvalues
$\lambda_2$ through $\lambda_{101}$, which remain $O(1)$, in the
Frobenius norm; only the ground state $\lambda_1$ decreases across
the 113-OOM range.
The bulk trace $\sum_{i=2}^{101}\lambda_i$ shows a slow monotonic
decay well-described by a power law $\operatorname{Tr}(Q)\propto c^{-0.116}$
($R^2=0.994$, measured across 15~cutoffs), a total drift of~$\sim$13\%
from $c=13$ to $c=67$.  The eigenvalue distribution shifts but the
$L^2$~norm of the operator is approximately preserved.  The spectral gap
$\lambda_2/\lambda_1\sim10^{7\text{--}8}$ indicates that the ground
state is extremely well-isolated from the bulk spectrum at all cutoffs.

Two further bulk invariants are tightly linear in $c$ and $L=\log c$
respectively.  The log-determinant
$\log\bigl|\det Q_c^{\mathrm{even}}\bigr| =
\sum_{k=1}^{101}\log|\lambda_k|$ is
\emph{linear in $c$}: across eight cutoffs with available full spectra
($c\in\{13,14,23,29,31,37,41,43\}$),
\begin{equation}\label{eq:logdet}
  \log\bigl|\det Q_c^{\mathrm{even}}\bigr|
  \;\approx\; -65.6\, c + 542,\qquad R^2 = 0.997.
\end{equation}
The slope $-65.6$ is a finite-matrix determinantal invariant: each
unit increment in $c$ removes approximately $65.6$ natural-log units
of measure from the product spectrum, dominated by the
exponential decay of $\lambda_{\min}$ relative to the bulk.

Complementarily, the heat-kernel trace
$\operatorname{Tr}\bigl(e^{-tQ_c^{\mathrm{even}}}\bigr) =
\sum_{k=1}^{101}e^{-t\lambda_k}$ is well-described by an affine function
of $L=\log c$ for fixed~$t$:
\begin{equation}\label{eq:heat-kernel}
  \log\operatorname{Tr}\bigl(e^{-tQ_c^{\mathrm{even}}}\bigr)
  \;\approx\; A(t)\,L + B(t),
\end{equation}
with linear-fit RMS residual $1.4\times 10^{-4}$ at $t=10^{-2}$ and
$1.2\times 10^{-3}$ at $t=10^{-1}$ across the same eight cutoffs.
The smoothness in $L$ is two orders of magnitude tighter than the
log-eigenvalue mean (whose smallest-eigenvalue contribution is itself
non-smooth in $c$); only the ground-state eigenvalue contributes
the non-smooth behavior, while the bulk integrand contributes a
controllable analytic function of $L$.

\subsection{Random-matrix statistics of the bulk spectrum}
\label{sec:rmt}

Montgomery's pair correlation conjecture \cite{Montgomery1973}
predicts GUE statistics for the non-trivial zeros of the Riemann zeta
function, but only after \emph{local rescaling}: for $\gamma_n$ at
height~$T$ one rescales
$\gamma_n \mapsto \gamma_n\,\log(\gamma_n/2\pi)/(2\pi)$
so that the mean spacing is unity, and the GUE signature appears in
the \emph{local} pair correlation of the rescaled zeros.  The
conjecture is not expected to manifest at the un-rescaled (infrared)
level, where the zero density itself grows like
$\log T$ and is highly non-uniform.  A direct numerical test of
Montgomery in this framework would therefore require computing
many~$\gamma_k$ zeros at a single large cutoff and rescaling them
locally before measuring spacings; that is beyond the scope of the
present sweep.

Instead, we examine here the \emph{bulk} eigenvalues of the CvS
Galerkin matrix (eigenvalues $\lambda_2$ through $\lambda_{101}$,
excluding the exponentially small ground state) at the un-rescaled
level, as a structural diagnostic of the truncated operator.  No
GUE behavior is expected at this level; we record the bulk-spectrum
statistics for structural interest only, not as a test of
Montgomery.

We analyze the even-sector Galerkin matrix (dimension~101) rather
than the odd sector or a smaller subset, to maximize the sample
size for spacing statistics.  After polynomial unfolding (degree~5)
to unit mean spacing, we compute the nearest-neighbor spacing
distribution (NNSD) and fit the Brody parameter $\beta$, which
interpolates Poisson ($\beta=0$, no level repulsion) and GOE
($\beta=1$, linear repulsion).

\begin{table}[ht]
\centering
\caption{Kolmogorov--Smirnov distance to theoretical distributions
  and Brody parameter for bulk eigenvalues ($\lambda_2$ through
  $\lambda_{101}$).}
\label{tab:rmt}
\begin{tabular}{@{}rcccc@{}}
\toprule
$c$ & $\beta_{\text{Brody}}$ & KS(Poisson) & KS(GOE) & KS(GUE) \\
\midrule
13 & $<0.05$ & 0.125 & 0.253 & 0.210 \\
23 & $<0.05$ & 0.108 & 0.174 & 0.232 \\
37 & $<0.05$ & 0.106 & 0.145 & 0.305 \\
43 & $<0.05$ & 0.108 & 0.126 & 0.326 \\
53 & $<0.05$ & 0.152 & 0.162 & 0.263 \\
\bottomrule
\end{tabular}
\end{table}

\begin{figure}[!htbp]
\centering
\includegraphics[width=0.85\textwidth]{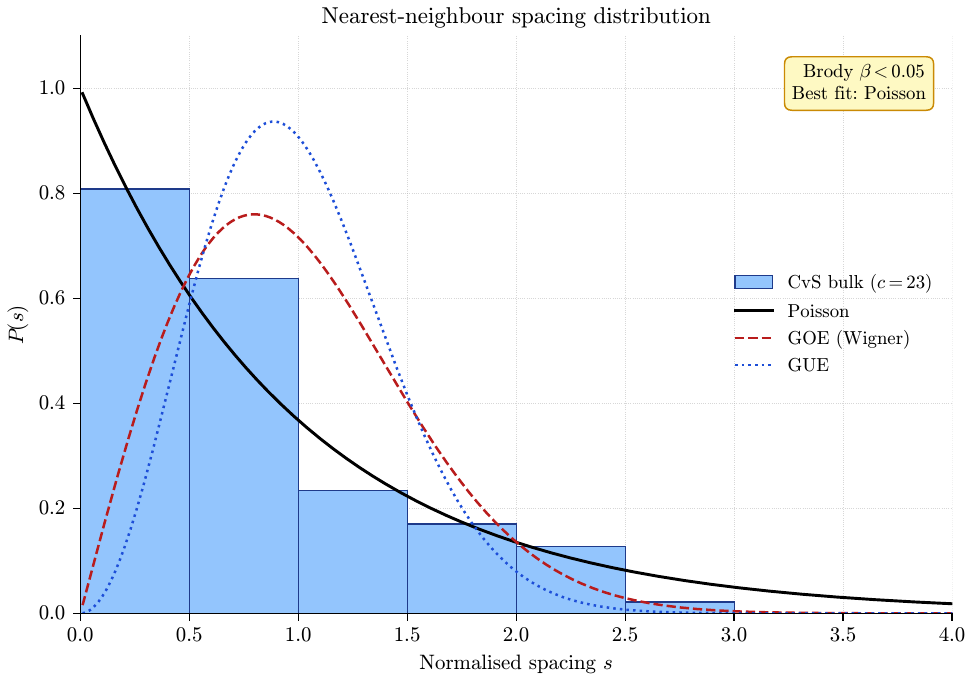}
\caption{Nearest-neighbor spacing distribution of the 100 bulk
  eigenvalues ($\lambda_2$ through $\lambda_{101}$) at $c=23$, after
  degree-5 polynomial unfolding to unit mean spacing.  The histogram
  matches the Poisson distribution (no level repulsion), consistent
  with Brody $\beta<0.05$ at all tested cutoffs.}
\label{fig:rmt-spacing}
\end{figure}

The result is uniform across the ten cutoffs tested ($c=13$ through
$c=53$; the $c=59$, $c=61$, and $c=67$ eigenvectors were computed
after this analysis): $\beta<0.05$ (consistent with Poisson,
negligible level repulsion) at every cutoff
(Figure~\ref{fig:rmt-spacing} shows a representative histogram at $c=23$).
KS tests indicate
that the Poisson distribution is the best fit at every cutoff,
with GOE second and GUE worst.

The Brody parameter is reported as $\beta<0.05$ at all cutoffs, likely
reflecting the limited sample size ($N_{\mathrm{bulk}}=100$ eigenvalues)
rather than exact Poisson statistics.

The Poisson classification indicates that the bulk eigenvalues behave
as uncorrelated levels --- characteristic of an integrable system.
This is physically consistent with the enormous spectral gap
$\lambda_2/\lambda_1\sim10^{7\text{--}8}$ isolating the ground state
from the bulk: the bulk eigenmodes are essentially independent of the
prime-sum perturbation that drives the ground-state convergence.

We emphasize three caveats: (i)~with $N=100$ bulk eigenvalues, these
are small-sample estimates; (ii)~the bulk eigenvalues are properties
of the \emph{finite-dimensional Galerkin truncation}, not of the
Riemann zeros themselves; and (iii)~per Montgomery, GUE statistics
for $\zeta$-zeros are a property of the \emph{locally-rescaled}
zeros and are not expected to manifest in the un-rescaled bulk
Galerkin spectrum studied here.  A proper random-matrix test of
Montgomery in this framework would require many extracted
$\gamma_k$ zeros at a single large cutoff, locally rescaled by the
mean zero density, which we do not undertake.

\subsection{Per-prime convergence decomposition}
\label{sec:per-prime}

The 15-point dataset permits a decomposition of the convergence gain
at each step into contributions from interval extension versus new
prime inclusion.  The $c=13\to14$ step is particularly informative:
$c=14$ is not prime, so $\pi(14)=\pi(13)=6$ and no new prime enters
the operator.  Yet the error improves by 5.8~OOM (from $10^{-54.70}$ to
$10^{-60.45}$), purely from the increase in interval length
$L=\log c$.  This step has the highest normalized efficiency
of the entire sweep.

For the 11~prime-adding steps ($c=14\to17$ through $c=61\to67$), the
per-step gain $\Delta\log_{10}|\gamma_1\text{ error}|$ varies by a
factor of~10 (from $-1.65$ to $-17.08$).  However, this variation is
almost entirely explained by the interval length change $\Delta L$:
\[
  \text{Pearson}\;r(\text{per-step gain},\;\Delta L) = -0.96
  \qquad (n=14,\;p<10^{-5}).
\]
(Some degree of correlation between step size and gain is expected for
any monotone process sampled at unequal intervals; the high $|r|$
reflects the dominance of this geometric effect over prime-specific
contributions.)
After normalizing by $\Delta L$, the efficiency
$\Delta\log_{10}|\text{error}|/\Delta\log_{10}c$ has coefficient of
variation 0.17, compared to 0.52 for the raw gains --- a threefold
reduction.  The efficiency declines monotonically, consistent with
diminishing marginal returns at higher cutoffs.

Neither the von~Mangoldt weight $\log p/\sqrt{p}$ nor the
prime-counting function $\pi(c)$ significantly predicts per-step gains
after controlling for $\Delta L$.  The convergence gain at each step
is driven by how much the interval $[0,L]$ grows, not by intrinsic
properties of the newly included prime.

\subsection{Step-size power law and asymptotic behavior}
\label{sec:step-power-law}

Analyzing the per-step convergence gains as a function of~$c$ reveals
a power-law decay in the marginal gain per unit of cutoff.  In the
late regime ($c\geq31$), the gain per step is well described by
\[
  |\Delta\log_{10}|\text{error}|| \approx 381\cdot\mathrm{gap}(c)\cdot c^{-1.48},
  \qquad R^2 = 0.996,
\]
where $\mathrm{gap}(c) = c - c_{\mathrm{prev}}$ is the gap to the
previous cutoff.

Since the exponent $\alpha=1.48>1$, the sum of all future step sizes
$\sum_{p>67}\mathrm{gap}(p)\cdot p^{-1.48}$ converges by the integral
test.  This raises the possibility that
$\log_{10}|\gamma_1\text{ error}|$ approaches a \emph{finite limit}
as $c\to\infty$.  Integrating the power law from $c=67$ onward gives
a rough estimate of $\sim\!106$ remaining orders of magnitude,
suggesting an asymptotic limit near $10^{-274}$.
This extrapolation rests on 8~data points in the late regime ($c\geq31$)
and should be treated as a speculative estimate, not a prediction.
The power law may break at higher cutoffs.

\textbf{Two interpretations.}
\begin{enumerate}
\item \textbf{Finite limit.}  If $\alpha$ remains above~1, the
  first-zero error converges to a nonzero value ($\sim\!10^{-274}$).
  The CvS operator would converge to a function whose first zero is
  extremely close to, but not exactly at, $\gamma_1$.
\item \textbf{Exponent drift.}  The exponent $\alpha=1.48$ is measured
  over only $c\in[31,67]$.  At larger cutoffs, $\alpha$ could decrease
  below~1, in which case the series diverges and
  $|\gamma_1\text{ error}|\to0$ (full convergence).  With 15~data points,
  these scenarios cannot be distinguished.
\end{enumerate}
Both interpretations are presented here for completeness.  The
distinction is directly testable by extending the computation to
$c\geq71$.

\subsection{\texorpdfstring{$\chi_3$}{chi\_3} odd-sector analysis}
\label{sec:chi3-odd}

At $c=23$ and $c=29$, computed at the archimedean cutoff $T=400$, the
even sector of the $\chi_3$ CvS operator shows a single negative
eigenvalue.  As at $c=100$ (Section~\ref{sec:c100-doublets}), this is a
finite-cutoff artifact: increasing the cutoff removes it (the
$-6.46\times10^{-23}$ at $c=23$ becomes positive by $T=800$, the
$-5.82\times10^{-17}$ at $c=29$ by $T=1200$), leaving the even sector
non-negative.  At $T=400$ the full 101-eigenvalue spectrum shows a
sharp sector dichotomy:

\begin{table}[ht]
\centering
\caption{Negative eigenvalue count in the $\chi_3$ CvS operator at
$T=400$.  The $c=23,29$ even-sector negatives are finite-cutoff
artifacts, absent at larger $T$.}
\label{tab:chi3-odd}
\begin{tabular}{@{}rcccl@{}}
\toprule
$c$ & \#negative & Negative eigenvalue & Odd sector & Status \\
\midrule
13 & 0 & --- & all positive & positive \\
17 & 0 & --- & all positive & positive \\
23 & 1 & $-6.46\times10^{-23}$ & all positive & even-sector only \\
29 & 1 & $-5.82\times10^{-17}$ & all positive & even-sector only \\
37 & 0 & --- & all positive & positive \\
\bottomrule
\end{tabular}
\end{table}

At every cutoff, the odd sector (50~eigenvalues) is strictly positive.
At $T=400$ this sign change is confined to a single even-sector
eigenvalue.  The negative eigenvalues at $c=23$ and $c=29$ are
17--23 orders of magnitude smaller than the typical bulk eigenvalue
scale ($\sim\!1$--7), far below the $T=400$ archimedean truncation
error in the matrix entries, which is why the finite cutoff tips their
sign; both become positive once $T$ is increased.

The genuine sector pattern is that the odd sector is strictly positive
at every cutoff tested.  The apparent even-sector positivity breakdown
at $c=23,29$ is a finite-cutoff artifact, not a character-dependent
feature of the CvS framework.

\section{Related Work}\label{sec:related}

\subsection{The Connes program}

The spectral interpretation of the Riemann zeros has been a central
theme in Connes' program since the late 1990s \cite{Connes1999}.  The
present work sits within the specific thread initiated by Connes'
2026 paper \cite{Connes2026}, which revisits the connection between
the Weil explicit formula and the spectral action principle.

\subsection{Connes--van Suijlekom (2025)}

\CvS\ \cite{CvS2025} provide the mathematical foundation: Theorem~6.1
establishes that the zeros of $\hat{\eta}_c$ lie on the critical line
for every finite~$c$.  The Galerkin matrix $q_{m,n}$ of
Proposition~4.1 is the object we implement.

\subsection{Connes--Consani--Moscovici (2025)}

\CCM\ \cite{CCM2025} present a rank-one perturbation of the scaling
spectral triple with Euler products over primes $p\leq\lambda^2=c$.
As established in Section~\ref{sec:setup}, the \CCM\ matrix
$\tau_{i,j}$ is unitarily equivalent to the \CvS\ matrix $q_{m,n}$.
\CCM\ Section~8 explicitly names the absence of a convergence rate as
the central open problem --- the same problem addressed by the
measurements in the present work.

\subsection{Connes--Consani (2020)}

\CC\ \cite{CC2020} provide the explicit formula for the archimedean
Mellin multiplier $h_+(\tau)$ used in our implementation.

\subsection{\CCM\ (2023)}

\CCM\ \cite{CCMprolate2023} introduce the prolate wave operator,
which appears in the Connes-program framework as the approximation
construction for the limit $k_\lambda$ (Connes 2026 \cite{Connes2026}
\S6.3--\S6.4 and the ``educated guess'' of \CCM\ \cite{CCM2025} \S7).
The \CCM\ 2025 Galerkin computations themselves use the trigonometric
basis, as is evident from \CCM\ \cite{CCM2025} Lemma~5.1.

\subsection{Montgomery (1973)}

Montgomery \cite{Montgomery1973} conjectured that the local pair
correlation of the non-trivial zeros of the Riemann zeta function,
\emph{after local rescaling by the mean zero density}, follows GUE
statistics.  Our bulk-spectrum analysis (Section~\ref{sec:rmt}) is
\emph{not} a test of Montgomery's conjecture: we record un-rescaled
finite-dimensional Galerkin bulk-spectrum statistics as a structural
diagnostic of the truncated operator, at a regime (the un-rescaled,
infrared level) where Montgomery's conjecture is not expected to
manifest.

\subsection{\'Sliwi\'nski (2026): a complementary numerical regime}

After our April 2026 literature sweep, \'Sliwi\'nski
\cite{Sliwinski2026} independently studied the spectra of the
\CCM\ rank-one operator $D_{\log}^{(\lambda,N)}$, deriving an
inverse-logarithmic lower bound on the dissonance between those
spectra and Riemann zeros in the joint regime $\lambda = N \to \infty$
(roughly $c = \lambda^2 = N^2$).  His implementation operates at
approximately seven-digit GPU precision and reaches $N$ up to several
thousand; his target statistic is an averaged dissonance metric over
many zeros, not the per-cutoff matching-digit accuracy on a fixed
$\gamma_k$ studied here.  The two works cover disjoint regimes
(joint $\lambda{=}N$ growth at low precision versus the
$N\to\infty$-first-then-$c$ Galerkin sweep at high precision) and
test different quantities.  At the time of the original submission we were
aware of no other independent public numerical study of the \CvS\ Galerkin
matrix in the trigonometric-basis matching-digit-precision regime explored
here; several have appeared since, and are recorded in
Section~\ref{sec:independent-repro}.

\subsection{Independent reproductions}\label{sec:independent-repro}

A literature sweep conducted in April 2026 and refreshed in May 2026
--- encompassing arXiv, Google Scholar, MathOverflow, GitHub, and
seminar video recordings --- found no independent public reproduction
of Connes' $c=13$ first-zero error at matching-digit precision beyond
\'Sliwi\'nski's complementary work just discussed.  The only
matching-digit corroborating data outside this project was from the
\CCM\ author group itself.  No public code implementing the
\CvS\ Proposition~4.1 Galerkin matrix at arbitrary precision was
found in any language.

\medskip
\noindent\textbf{Note added (2026 August).}  That situation has changed since
the first version of this paper.  Five parties have since reported
reproductions: three from implementations independent of the code accompanying
this paper, one a cross-platform rerun of that package, and one a
reproduction of a single published table row with an independent harness.
Each item states what its author reports; the listing is attribution rather
than endorsement, and none of these results has been re-derived here, the one
exception being noted in the fourth item.
\begin{itemize}
\item B.~Martin (Skyline Trail Computing) reimplemented the \CvS/\CCM\
  Galerkin matrix from scratch, as a separate multiprecision assembly sharing
  no code with this project.  At the $c=100$, $N=250$, $\dps=500$ extraction
  cell he reports ordinates agreeing with those detected here from roughly
  $330$ digits at $\gamma_1$ to $\sim309$ at $\gamma_{10}$; at the $c=13$,
  $N=80$, $T=400$, $\dps=80$ cell he reports recovering $\gamma_1$ with
  absolute error $1.77\times10^{-55}$ against the true ordinate, and reports
  his separate comparison against the published \CCM/Connes error values at
  that cutoff as agreeing to a factor of $1.2$--$1.3$.
\item R.~Andrews reports an independent reproduction of the $c=13$ and
  $c=100$ spectra, together with a convergence analysis (Zenodo concept DOI
  \texttt{10.5281/zenodo.20427499}, at the time of writing resolving to
  record \texttt{21725468}).
\item A.~F.~Martini reports an independent from-scratch \texttt{mpmath}
  implementation of the \CCM\ construction (Zenodo
  \texttt{10.5281/zenodo.21864192}): strict positivity of the smallest
  eigenvalue at $c=13,17,19,23,29$ at $N=60$, taken at $\dps=80$ for
  $c=13,17,19$ and at $\dps=100$ for $c=23,29$, where he reports that
  $\dps=80$ returned negative values below its own resolution floor and
  diagnoses them as precision artifacts; recovery of the first ten ordinates
  at $c=13$ with errors of $1.4$--$5.0\times10^{-17}$, a level he attributes
  to finite truncation at $N=60$ rather than to the cutoff; and a check that the
  \CCM\ three-term decomposition and the \CvS\ Cauchy--Toeplitz assembly are
  the same form in different bases, to $2\times10^{-42}$.
\item The GitHub user \texttt{karl-keysingularity} reports rerunning the
  accompanying package (version 0.2.2) at the $c=13$, $N=100$, $T=400$,
  $\dps=80$ cell under native Windows and deposited the resulting values as a
  machine-readable artifact.  This is a cross-platform reproduction of the
  package rather than an independent implementation.  Comparing his deposited
  values against the committed reference gives $22$ significant digits on
  $\lambda_{\min}^{\mathrm{even}}$ and all $20$ stored digits of the
  $\gamma_1$ error; those two digit counts were computed here from the values
  he supplied.
\item M.~Osman reports agreement in every printed digit of the $c=23$ row of
  Table~\ref{tab:sobolev-scaling}, measured with his own harness, in the
  course of the three-cutoff test of Section~\ref{sec:sobolev-pw}.  This is a
  printed-decimal agreement claim on one row, not a spectrum reproduction.
\end{itemize}
\noindent The sweep result above should accordingly be read as a record of
the literature as it stood in mid-2026, not as a current claim.  Public code
for the \CvS\ Galerkin matrix now exists: it is the package accompanying this
paper (Section~\ref{sec:code}).

\section{Limitations and Honest Non-Claims}\label{sec:limitations}

We list the limitations of this work explicitly:
\begin{itemize}
\item \textbf{We do not claim to validate Connes' conjecture.}  The
  agreement at $c=13$ is a consistency check, not a comparison against
  an independent prediction.

\item \textbf{We do not claim to measure a convergence rate.}  All eight
  two-parameter smooth-convergence fits failed at 15~data points, and
  a five-parameter log-periodic model was falsified out-of-sample.
  Whether this reflects genuine non-smoothness or the need for
  different model families is unresolved.

\item \textbf{We do not claim any progress on RH.}  The \CvS\
  framework's Hurwitz-type convergence argument is a conjecture.  We
  compute some numbers in the neighborhood of that conjecture without
  bearing on its truth.

\item \textbf{The factor-of-1.3 discrepancy at $c=13$ is unresolved.}
  Plausible sources include precision, $N$, normalization conventions,
  integration cutoff~$T$, and root-finding tolerance.

\item \textbf{We compute the even sector only.}  The odd-sector zero,
  if any, is not probed here.  It has since been probed independently by
  M.~Osman, whose sign-change experiments on the odd-sector ground
  eigenvector we have reproduced.

\item \textbf{The $N=100$ multi-cutoff sweep is at $c\leq67$.}  The
  cross-extension --- an $N=100$ datum at $c\geq 71$ that would
  separate the $c$-dependent and $N$-dependent finite-$N$ effects
  more cleanly at the lower end of the sweep --- is left to future
  work.  The $c=100$ measurements include $N\in\{100,150,200,250\}$
  at $\dps=500$, providing four points on the $N$-sweep, sufficient
  to test the geometric-convergence hypothesis on consecutive Aitken
  triples (Section~\ref{sec:c100-aitken}).

\item \textbf{The Galerkin truncation bias at fixed $N$ is not
  independently measurable.}  We cannot distinguish the ``true
  eigenvalue at infinite~$N$'' from our $N=100$ eigenvalue without
  either a theoretical bound or a prolate-basis implementation (which,
  as noted in Section~\ref{sec:related}, would be a genuinely distinct
  direction not yet reported by any group).

\item \textbf{The pre-registered five-bin classifier fired ``Ambiguous
  partial floor.''}  This is an intermediate verdict.

\item \textbf{The step-size power law ($\alpha=1.48$) is measured over
  a limited range.}  Whether $\alpha$ remains above~1 at larger cutoffs
  is unknown; the finite-limit interpretation is one of two
  possibilities.

\item \textbf{CCM Theorem 1.1 is applied to an empirically
  distinguished positive branch at $c=100$, not to the raw matrix
  spectrum.}  At $c \leq 67$, $N=100$ the finite-$N$ Galerkin matrix
  is observed positive on every cell with the smallest eigenvalue
  numerically simple and an even-dominant eigenvector, consistent
  with the \CCM\ Theorem~1.1 ``$\epsilon_N$ simple, $\xi$ even''
  hypotheses as numerical conditions (not as analytically
  established hypotheses of the continuum theorem).  At $c=100$ with $N \geq 100$, the raw Galerkin
  spectrum computed at the finite archimedean cutoff $T=800$ contains
  negative-sign eigenvalues (Section~\ref{sec:c100-doublets}); these
  are an artifact of that cutoff, absent once $T$ is increased, so
  cutoff-free the $c=100$ even sector is non-negative and its smallest
  eigenvalue coincides with the smallest-\emph{positive} branch whose
  eigenvector gives the reported zero extraction
  (Table~\ref{tab:c100-gamma}).  The $c=100$ data are therefore
  consistent with the \CCM\ Theorem~1.1 ground-state hypotheses on the
  same footing as the $c\leq67$ cells.  Continuum positivity
  of $QW_\lambda$ is RH-equivalent and is not assumed at $\lambda
  = \sqrt{100}$.  See Remark~\ref{rem:ccm-hyp-c100}.  The $c=67$,
  $N=150$ corroborative measurement remains in the matrix-positive
  regime where the formal hypothesis still holds.

\item \textbf{The independence claim requires calibration.}  This
  implementation was drafted with AI assistance under continuous
  human supervision (see the Statement on use of AI tools, p.~\pageref{sec:ai-tools}).
  It is independent of the Connes/\CCM\ codebase in the sense that
  no part was derived from non-public code.  It is not independent in the
  sense of having been written by a second human research group.  Two
  from-scratch reimplementations by other parties have since been reported
  and are recorded in Section~\ref{sec:independent-repro}; as reported, they
  address this limitation from the outside rather than removing it here.

\item \textbf{No part of this paper constitutes a proof of anything.}
  This is experimental mathematics in the literal sense.
\end{itemize}

\section{Future Directions}\label{sec:future}

\paragraph{Extension to Dirichlet $L$-functions.}
We extended the CvS pipeline to $L(s,\chi_3)$, where $\chi_3$ is the
non-principal character modulo~3, at five cutoffs ($c=13,17,23,29,37$).
The archimedean piece shifts from $\psi(\tfrac14+i\tau/2)$ to
$\psi(\tfrac34+i\tau/2)$ for this odd character, and the prime-piece
weights are multiplied by $\chi_3(n)$.

The results show a non-monotonic first-zero error:

\begin{table}[ht]
\centering
\caption{$L(s,\chi_3)$: first-zero error across cutoffs.  The
$\lambda_{\min}$ signs are computed at $T=400$; the negatives at
$c=23,29$ are finite-cutoff artifacts (Section~\ref{sec:chi3-odd}).}
\label{tab:chi3}
\begin{tabular}{@{}rccr@{}}
\toprule
$c$ & $\lambda_{\min}$ sign & $|\gamma_1(\chi_3)\text{ err}|$ & digits \\
\midrule
13 & $+$ & $4.18\times10^{-17}$ & 16 \\
17 & $+$ & $7.09\times10^{-25}$ & 24 \\
23 & $-$ & $2.00\times10^{-23}$ & 22 \\
29 & $-$ & $5.37\times10^{-18}$ & 17 \\
37 & $+$ & $7.30\times10^{-29}$ & 28 \\
\bottomrule
\end{tabular}
\end{table}

At $c=13,17,37$ the even-sector minimum eigenvalue is positive; the
first-zero error is non-monotonic in $c$ (16, 24, 22, 17, 28~digits at
$c=13,17,23,29,37$).  The $\lambda_{\min}$ sign shown is the $T=400$
value: the negatives at $c=23,29$ are a finite-cutoff artifact
(Section~\ref{sec:chi3-odd}), positive once $T$ is increased, the same
archimedean-truncation mechanism that produces the $c=100$ block of
Section~\ref{sec:c100-doublets}.  Whether the modest convergence dips
at $c=23,29$ reflect a genuine character effect or the finite-$T$ entry
error is not settled here, and is a natural question for extending the
CvS framework to general Dirichlet $L$-functions.

\paragraph{Predictions for future cutoffs.}
At $\dps=200$, the backward-error floor permits clean measurements
to approximately $c=97$ before $\dps=300$ would be required.  We
predict continued monotone improvement in $|\gamma_1\text{ error}|$
with each new cutoff, with the $|\gamma_1\text{ err}|/\lambda_{\min}$ ratio
continuing its slow upward drift.  The $N$-convergence data suggest
that $N=100$ is sufficient for cutoffs up to at least $c=53$; whether
larger~$N$ is needed at $c\geq59$ is an open empirical question.

\paragraph{Multi-zero convergence at higher cutoffs.}
The multi-zero universality established in
Section~\ref{sec:multi-zero} (all ten zeros within 3.8\% of the
$\gamma_1$ rate) should be tested at cutoffs beyond $c=67$ and with
larger basis sizes ($N>100$) to determine whether the slight
systematic slowdown ($\sim\!0.42\%$ per zero index) persists or
saturates.

\paragraph{Transition from Poisson to GUE?}
The uniformly Poisson bulk spectrum observed at $N=100$ raises the
question of whether GUE statistics emerge at larger~$N$ or higher
cutoffs.  Montgomery's conjecture predicts GUE for the Riemann zeros
themselves; whether the finite-dimensional CvS Galerkin truncation
approaches GUE statistics in any limit is unknown.

\paragraph{Theoretical bound on $C(c)$.}
The ratio $|\gamma_1\text{ err}|/\lambda_{\min}\approx C(c)$ with
$C(c)\approx6730\log c - 11268$ is an empirical fit.  A theoretical
derivation of the growth rate of $C(c)$ from the structure of the
\CvS\ operator would constitute a substantive step in understanding
the convergence mechanism.

\paragraph{Testing the Sobolev--Paley--Wiener prediction: done, and
negative.}
This section previously pre-registered the following test, quoted from the
earlier version: ``The testable prediction of Section~\ref{sec:sobolev-pw},
that the Sobolev slope should scale linearly with~$T$, can be verified by
repeating the $N$-convergence measurement at $T=1600$.  If the slope doubles
from ${\sim}55$ to ${\sim}110$, the Paley--Wiener mechanism is supported; if
it remains constant, the mechanism is incorrect.''  The test has been run at
$c=23$ on the $N\in\{40,60,80\}$ grid at $\dps=150$ and the exponent is
constant to within $0.206$ across $T=400,800,1600$
(Section~\ref{sec:sobolev-pw}).  The mechanism is accordingly withdrawn.  A
single cutoff suffices for that conclusion, since proportionality is asserted
at fixed~$c$; what remains open is a correct account of the observed
$c$-dependence of the exponent, for which we have no candidate mechanism.

\paragraph{Resolving the step-size exponent.}
Computing $c=71,73,79$ would test whether the step-size power-law
exponent $\alpha\approx1.48$ remains above~1, which would support the
finite-limit interpretation, or drifts below~1, which would support
full convergence to zero error.

\section{Code and Data Availability}\label{sec:code}

A complete reference implementation is available as the open-source
Python package \texttt{connes-cvs} (version 0.2.2 at the time of the
original submission), installable via \texttt{pip install connes-cvs} from
PyPI.  Source: \url{https://github.com/akivag613/connes-cvs-}.
The ancillary files attached to this arXiv submission contain the
principal result JSONs and the driver scripts that reproduce each
table and figure in this paper.

\paragraph{Package version note.}
The current release of \texttt{connes-cvs} is the one named in the
package \texttt{CHANGELOG} at the repository above.  The 0.3 series
changes no value reported here.  Two of its changes are worth stating
for anyone reproducing
from the paper.  First, the zero-extraction entry point now accepts the
cutoff~$c$ directly and forms $L=\log c$ internally at the working precision;
passing~$L$ as a double-precision float is still accepted but now warns,
because a float64 $L$ silently caps extraction accuracy near $10^{-16}$
irrespective of the working precision.  That defect was reported independently by
M.~Osman and by the GitHub user \texttt{karl-keysingularity}.  In 0.2.2 it
affected the quick-start example and the $\gamma_k$ errors reported by
\texttt{connes\_cvs.sweep.run\_sweep}, which the ancillary drivers
\path{A_extended_c_sweep.py}, \path{U_T800_full_sweep.py} and
\path{X_c67_blind_test.py} wrap; it did not touch matrix assembly or the
eigensolve.  It did not affect the values reported in this paper, which were
produced by standalone full-precision drivers and whose errors run from
$10^{-55}$ to $10^{-168}$, far below the $10^{-16}$ cap the defect imposes.
Anyone re-running the ancillary scripts should use 0.3.0 or later; under 0.2.2 they
return a $\gamma_1$ error capped near $10^{-16}$ irrespective of the working
precision.  Second, the supported Python floor rises from 3.9 to 3.10.  Version 0.2.2 remains installable for
byte-level comparison against the originally published cells.  Corrections to
this paper are tracked in \path{ERRATA.md} in the repository above and in the
Zenodo record.

\paragraph{Ancillary files.}
All files are deposited at the top level of the Zenodo record (flat
layout, individually previewable).  The \LaTeX{} source of this
manuscript is distributed by the arXiv submission, not by this Zenodo
deposit; figures are embedded in the PDF above and need not be
regenerated for reproduction of the numerical claims.

\emph{Numerical results.}
\begin{sloppypar}\raggedright
\begin{itemize}
\item \path{results_15pt_T800.json} --- the
  15-cutoff sweep (Table~\ref{tab:15pt-sweep}), $c=13$ through $67$ at
  $N=100$, $T=800$.
\item \path{results_N_sweep_c23.json} ---
  $N$-convergence at $c=23$ (Table~\ref{tab:N-conv}).

\item \path{results_c100_N100_T800_dps500_v020.json},
  \path{results_c100_N150_T800_dps500_v020.json},
  \path{results_c100_N200_T800_dps500_v020.json},
  \path{results_c100_N250_T800_dps500_v020.json},
  \path{results_c100_N150_T800_dps1000_v020.json}
  --- the four-point $N$-sweep at $c=100$ plus the $\dps=1000$ retest
  at $N=150$ (Table~\ref{tab:c100-N-sweep}).
\item \path{c100_N150_dps500_gamma_extraction.json},
  \path{c100_N150_dps1000_gamma_extraction.json},
  \path{c100_N250_dps500_gamma_extraction.json}
  --- $\gamma_1,\ldots,\gamma_{10}$ extractions for the three precision
  cells of Table~\ref{tab:c100-gamma}; the $N=250$ file stores both
  detected and reference values to 400 significant digits for
  independent verification.
\item \path{c100_full_spectrum_diagnosis.json}
  --- all 151 eigenvalues at $c=100$, $N=150$, $\dps=500$, with signs;
  source for the finite-$T$ ($T=800$) negative-sign artifact block of
  Section~\ref{sec:c100-doublets}.
\item \path{power_law_refit.json} --- the
  $13.24\,c^{0.634}$ in-sample fit on $c\leq67$, the $c=100$ residual
  of $49$~OOM that falsifies it, and the four-parameter refit including
  $c=100$ (Section~\ref{sec:c100-reframe}).
\item \path{richardson_n_extrapolation.json} ---
  legacy 3-point Aitken/Richardson analysis at $N\in\{100,150,200\}$;
  the 4-point analysis is recomputed by
  \path{c100_aitken_check.py}.
\item \path{results_c67_N150_T800_dps500_v020.json}
  --- the $c=67$, $N=150$ corroborative datum of
  Section~\ref{sec:c100-reframe}.

\item \path{eigenvector_overlaps.json} --- $15\times15$
  pairwise eigenvector overlap matrix (Section~\ref{sec:structural}).
\item \path{results_rmt_analysis.json} --- bulk-spectrum
  random-matrix statistics (Section~\ref{sec:rmt}).
\item \path{results_L_multi_zero.json} --- multi-zero
  convergence data across cutoffs (Section~\ref{sec:multi-zero},
  Figure~\ref{fig:multi-zero}).
\item \path{results_A_c23.pickle} --- full $c=23$ even-sector
  eigenspectrum and supporting metadata; source for the nearest-neighbour
  spacing distribution of Figure~\ref{fig:rmt-spacing}.

\item \path{preregistered_predictions_2026-04-12.md}
  --- the Section~\ref{sec:blind-c67} blind-prediction document, locked
  by SHA-256 hash before the $c=67$ computation was performed.
\end{itemize}
\end{sloppypar}

\emph{Driver scripts.}
\begin{sloppypar}\raggedright
\begin{itemize}
\item \path{c100_experiment_optimized.py} --- $c=100$ cell-based
  runner for the five cells in Table~\ref{tab:c100-N-sweep}
  (\path{c100-N100}, \path{c100-N150-dps500}, \path{c100-N200},
  \path{c100-N250-dps500}, \path{c100-N150-dps1000}).  Uses the
  $v0.2.0$ fused-kernel math core from the public \texttt{connes-cvs}
  PyPI package with \texttt{multiprocessing.imap\_unordered} for the
  $\psi$-cache phase.
\item \path{A_extended_c_sweep.py}, \path{U_T800_full_sweep.py}
  --- 15-cutoff $c$-sweep at $T=800$ (Section~\ref{sec:extension});
  thin wrappers around \texttt{connes\_cvs.sweep.run\_sweep}.
\item \path{X_c67_blind_test.py} --- single-cell $c=67$ runner
  (Section~\ref{sec:blind-c67} blind-prediction cell at $N=100$, $\dps=200$;
  Section~\ref{sec:c100-reframe} corroborative cell at $N=150$, $\dps=500$).
\item \path{c100_aitken_check.py} --- ${<}\,1$-second headline
  reproducer.  Loads the four $c=100$ JSONs and recomputes the
  Aitken-$\Delta^2$ extrapolation on both overlapping triples, prints
  the under-$1\%$-of-exponent agreement on the deeper triple with
  the Connes 2026~\S6.4 heuristic continuum prediction.
\end{itemize}
\end{sloppypar}

\paragraph{Package requirements} (from \texttt{pyproject.toml}):
\texttt{mpmath}~$\geq 1.3$, \texttt{python-flint}~$\geq 0.5$,
\texttt{gmpy2}; Python~3.10+.  The specific environment used for
the runs reported here was \texttt{mpmath}~1.4.1,
\texttt{python-flint}~0.8.0 on macOS with a 12-core Apple~M2~Max.
The \texttt{python-flint} digamma primitive provides an
$\sim 11\times$ speedup over pure mpmath for the archimedean
quadrature dominant in single-cell runs.

\paragraph{Reproducing the headline results.}
\begin{verbatim}
pip install connes-cvs
python examples/basic_compute.py          # c=13, ~5 min with python-flint
python examples/c100_aitken_check.py      # c=100 Aitken match, < 1 s
\end{verbatim}
The first command produces $|\gamma_1\text{ error}|\approx 1.5\times
10^{-55}$ at $c=13$; the second loads the ancillary $c=100$ data and
recomputes the Aitken-$\Delta^2$ extrapolation and the Connes 2026
\S6.4 prediction, printing the under-$1\%$-of-exponent agreement
on the deeper-anchored triple reported in
Section~\ref{sec:c100-aitken}.

\paragraph{Reproducibility certification.}  The example scripts above
use $T=400$, $\dps=80$ for speed (single-process wall-clock ${\sim}5$
minutes on Apple~M-series with python-flint installed), and produce
$|\gamma_1\,\mathrm{err}|\approx 1.455\times10^{-55}$ at $c=13$
(Table~\ref{tab:T-convergence}, $T=400$ row).  The paper's headline $c=13$
value $|\gamma_1\,\mathrm{err}|=2.005\times10^{-55}$ (Table~\ref{tab:15pt-sweep})
is the $T=800$, $\dps=200$ datum; the shift from $T=400$ to $T=800$ is
$0.14$ units in $\log_{10}|\gamma_1\,\mathrm{err}|$ (Table~\ref{tab:T-convergence}),
within the factor-of-1.3 envelope of the \CvS/\CCM\ measurements.
The $c=100$ cell-based runner
\path{c100_experiment_optimized.py} produces a JSON file whose
\path{lambda_even} field, taken as a decimal string and passed
through $\log_{10}(\cdot)$, reproduces the rows of
Table~\ref{tab:c100-N-sweep} bit-identically to working precision.
The $\gamma_k$-extraction JSON at $N=250$, $\dps=500$ stores both the
detected $\gamma_k$ and the reference \texttt{mpmath.zetazero($k$)}
value to $400$ significant digits, sufficient to independently
recompute $\lfloor-\log_{10}|\gamma_k-\gamma_k^{\mathrm{exact}}|\rfloor$
and verify the $307$--$329$ matching-digit counts of
Table~\ref{tab:c100-gamma} without rerunning the diagonalization.

\section*{Statement on use of AI tools}\label{sec:ai-tools}

In accordance with the arXiv policy on authors' use of generative AI
language tools, the author reports that Anthropic's Claude was used as
a text-to-text assistant for the drafting of \LaTeX{} exposition and
Python implementation code, to specifications written by the author.
No mathematical proof, theorem, or definition appearing in this paper
was originated by the language model; load-bearing mathematical claims
were verified by the author before inclusion.  The author conceived,
directed, and verified all computational experiments and mathematical
claims, and takes full responsibility for the contents of this
manuscript.

\section{Acknowledgments}\label{sec:ack}

I am grateful to Professor Alain Connes for substantive correspondence
in 2026 that confirmed the trigonometric-basis convention used in the
Connes--Consani--Moscovici Galerkin computations cited throughout this
paper, identified \CC\ \cite{ConnesConsani2023} \S3 as the qualitative
motivation underlying the $k_\lambda$ construction of Connes
2026~\cite{Connes2026} \S6.4 and its open-status approximation step
in \S6.6, offered to assist with the arXiv submission, and provided
specific formulation comments on the manuscript that improved the
abstract, the discussion of the $c=13$ cross-validation in
Section~\ref{sec:c13-1p3-discrepancy}, and the framing of the
bulk-spectrum analysis in Section~\ref{sec:rmt}.  The manuscript
benefited materially from his comments; the contents and any
remaining errors are the author's own.  I further thank Professor
Connes, whose questions about the $c=100$ spectrum prompted the
precision investigation behind the correction recorded in the note
added in revision; B.~W.~A.~Silva (Zenodo
\texttt{10.5281/zenodo.20650146}), who independently identified the
finite-cutoff sensitivity of the deep spectrum; and R.~Andrews (Zenodo
\texttt{10.5281/zenodo.20427500}), whose computation exhibits the
naturally even, positive ground state consistent with that finding.

I thank M.~Osman, who executed the test pre-registered in
Section~\ref{sec:future} and reported the result that refutes the
Paley--Wiener mechanism of Section~\ref{sec:sobolev-pw}, and who also
identified the cutoff mislabelling of Table~\ref{tab:sobolev-scaling}.  I
thank B.~Martin, R.~Andrews, A.~F.~Martini, M.~Osman and the GitHub user
\texttt{karl-keysingularity} for the reproductions recorded in
Section~\ref{sec:independent-repro}, and \texttt{karl-keysingularity} and
M.~Osman additionally for the independent reports of the float64 cutoff defect
that led to the precision-safety changes in release~0.3.0 of the accompanying
package.

The mathematical foundation of this work is entirely due to Connes,
van~Suijlekom, Consani, and Moscovici, through the papers cited in
Section~\ref{sec:related}.  No part of the implementation was derived
from non-public code held by any of these authors, and no public code
for the \CvS\ Galerkin matrix was located in the literature sweep
described in Section~\ref{sec:related}.

I am also grateful for substantive correspondence on the underlying
mathematics or on the code with Dr.\ Pascal~Molin (Universit\'e
Paris~Cit\'e, author of \texttt{Arb}'s \texttt{acb\_dirichlet}
primitive), Dr.\ Teun~van~Nuland (TU~Delft), and Dr.\ Sabine~Boegli
(Durham), whose readings prompted clarifications throughout.  I am
separately grateful to Professor Nigel~Higson (Penn~State) for the
warm referral that opened the line of correspondence with Professor
Connes, and to Professor Klaas~Landsman (Radboud) for a kind referral
when the topic fell outside his own area.

\bibliographystyle{plain}
\bibliography{paper1_v3}

\clearpage

\appendix

\section{Summary of Numerical Results}\label{app:results}

For reference, we collect the final precision-resolved measurements in
a single table.  All rows use $T=800$, $N=100$,
$\dps=150$ ($c\leq37$) or $\dps=200$ ($c\geq41$).

\begin{table}[ht]
\centering
\caption{Complete 15-point results at $N=100$, $T=800$,
  $\dps=150$ ($c\leq37$) / $\dps=200$ ($c\geq41$).}
\label{tab:appendix-results}
\begin{tabular}{@{}rllllr@{}}
\toprule
$c$ & $L=\log c$ & $\lambda_{\min}^{\mathrm{even}}$
    & $|\gamma_1\text{ error}|$
    & $|\gamma_1\text{ err}|/\lambda_{\min}$
    & $\log_{10}|\gamma_1\text{ err}|$ \\
\midrule
13 & 2.565 & $2.865\times10^{-59}$
   & $2.005\times10^{-55}$  & 6999  & $-54.70$ \\
14 & 2.639 & $4.835\times10^{-65}$
   & $3.541\times10^{-61}$  & 7324  & $-60.45$ \\
17 & 2.833 & $2.030\times10^{-80}$
   & $1.634\times10^{-76}$  & 8047  & $-75.79$ \\
19 & 2.944 & $1.265\times10^{-90}$
   & $1.070\times10^{-86}$  & 8457  & $-85.97$ \\
23 & 3.135 & $5.959\times10^{-107}$
   & $5.520\times10^{-103}$ & 9262  & $-102.26$ \\
29 & 3.367 & $4.366\times10^{-124}$
   & $4.587\times10^{-120}$ & 10507 & $-119.34$ \\
31 & 3.434 & $1.045\times10^{-128}$
   & $1.141\times10^{-124}$ & 10919 & $-123.94$ \\
37 & 3.611 & $4.670\times10^{-140}$
   & $5.686\times10^{-136}$ & 12177 & $-135.25$ \\
41 & 3.714 & $2.122\times10^{-146}$
   & $2.760\times10^{-142}$ & 13004 & $-141.56$ \\
43 & 3.761 & $2.519\times10^{-149}$
   & $3.379\times10^{-145}$ & 13412 & $-144.47$ \\
47 & 3.850 & $2.994\times10^{-154}$
   & $4.270\times10^{-150}$ & 14260 & $-149.37$ \\
53 & 3.970 & $9.615\times10^{-161}$
   & $1.493\times10^{-156}$ & 15529 & $-155.83$ \\
59 & 4.078 & $2.328\times10^{-166}$
   & $3.911\times10^{-162}$ & 16800 & $-161.41$ \\
61 & 4.111 & $5.063\times10^{-168}$
   & $8.722\times10^{-164}$ & 17226 & $-163.06$ \\
67 & 4.205 & $7.993\times10^{-173}$
   & $1.478\times10^{-168}$ & 18489 & $-167.83$ \\
\bottomrule
\end{tabular}
\end{table}

\begin{table}[ht]
\centering
\caption{Log-scale steps between consecutive cutoffs ($T=800$).}
\label{tab:steps}
\begin{tabular}{@{}lrr@{}}
\toprule
Pair & $\Delta\log_{10}|\gamma_1\text{ error}|$ & New prime \\
\midrule
$c=13\to14$ & $-5.75$ & (none) \\
$c=14\to17$ & $-15.34$ & 17 \\
$c=17\to19$ & $-10.18$ & 19 \\
$c=19\to23$ & $-16.29$ & 23 \\
$c=23\to29$ & $-17.08$ & 29 \\
$c=29\to31$ & $-4.60$ & 31 \\
$c=31\to37$ & $-11.31$ & 37 \\
$c=37\to41$ & $-6.31$ & 41 \\
$c=41\to43$ & $-2.91$ & 43 \\
$c=43\to47$ & $-4.90$ & 47 \\
$c=47\to53$ & $-6.46$ & 53 \\
$c=53\to59$ & $-5.58$ & 59 \\
$c=59\to61$ & $-1.65$ & 61 \\
$c=61\to67$ & $-4.77$ & 67 \\
\bottomrule
\end{tabular}
\end{table}

The step sizes are non-uniform, ranging from $-1.65$ ($c=59\to61$) to
$-17.08$ ($c=23\to29$).  The $c=13\to14$ step is notable: $c=14$ is not
prime, so no new prime enters the operator, yet the error improves by
5.8~OOM purely from the increase in interval length $L=\log c$.
The variation in step sizes is strongly correlated with $\Delta L=\Delta\log c$
(Pearson $r=-0.96$).  The normalized slope
$\Delta\log_{10}|\text{error}|/\Delta L$ is far more uniform than the
raw steps (coefficient of variation 0.17 vs~0.52).

\section{Pre-Registered Convergence Model Residuals}
\label{app:models}

The original three models (M1--M3) were pre-registered and fit to
3-point $\dps=80$ data in Iteration~5; all failed with max residuals
2.44--2.68.  Two additional models (M4, M5) were added for the
10-point sweep; three prime-counting models (M6--M8) were added for
the 15-point dataset.  A five-parameter log-periodic model was fit to
the 10-point training data and tested out-of-sample at $c=47$ and $c=53$.

\begin{table}[ht]
\centering
\caption{Full model-fit comparison: 15-point, 10-point, and 3-point
  residuals.  All models use 2~free parameters except the log-periodic
  model (5~parameters).}
\label{tab:model-comparison}
{\small
\begin{tabular}{@{}llccc@{}}
\toprule
Model & Form & Max res.\ (15-pt) & Max res.\ (10-pt)
      & Max res.\ (3-pt) \\
\midrule
M1 & $a\cdot L+b$             & 2.66  & 2.66 & 2.44 \\
M2 & $a\cdot L^2+b$           & 4.71  & 4.71 & 2.68 \\
M3 & $a\cdot L\log L+b$       & 3.63  & 3.63 & 2.56 \\
M4 & $a\cdot L+b\cdot\log L$  & 2.66  & 2.66 & --- \\
M5 & $a\cdot\exp(b\cdot L)$   & 1.77  & 1.77 & --- \\
M6 & $a\cdot\pi(c)+b$         & 10.28 & ---  & --- \\
M7 & $a\cdot\vartheta(c)+b$   & 12.09 & ---  & --- \\
M8 & $a\cdot\#\mathrm{pp}(c)+b$ & 10.20 & --- & --- \\
\midrule
\multicolumn{5}{@{}l}{\textit{Higher-parameter models:}} \\
Log-periodic (5-param)
   & $a\cdot L+b+A\sin(\omega L+\phi)$
   & \textbf{5.71} & 0.74 & --- \\
\bottomrule
\end{tabular}
}
\end{table}

M5 (exponential-in-$L$) is the best-performing two-parameter model
but still exceeds the $\leq0.5$ threshold by a factor of~3.5.
M6--M8, which use prime-counting independent variables, perform
substantially worse (max residuals 10--12), indicating that the
convergence is not proportional to any standard prime-counting
function.

The log-periodic model (5~parameters on 10~training points) achieved
max residual 0.74 on training data --- the first model to beat the
$\leq1.77$ threshold.  However, pre-registered blind predictions at
$c=47$ and $c=53$ yielded residuals of 1.5 and 5.5 respectively: the
model predicted that $c=53$ would be \emph{less} precise than $c=47$,
contradicting the observed monotone improvement.  The model was
overfitting noise in the training range.  This clean out-of-sample
falsification strengthens the negative result: no parametric model
tested --- including models with up to 5~free parameters --- captures
the $c$-dependence of $|\gamma_1\text{ error}|$.

\section{Reproducibility Discipline}\label{app:reproducibility-discipline}

This appendix collects the procedural rules enforced across the
project's iterations.  These are reproducibility-hygiene practices,
not mathematical claims; they are recorded here for transparency
about how the numerical results in the body were produced and
verified.

\begin{itemize}
\item \textbf{Pre-registration.}  Each iteration's experimental plan,
  including success/failure criteria and gate rules, is written and
  frozen before code execution.
\item \textbf{Regression checking.}  The $c=13$ row is verified
  to all reported digits against the prior iteration's result at
  each new iteration.
\item \textbf{Banned-phrase discipline.}  A cumulative list of banned
  phrases (including ``validates,'' ``confirms conjecture,'' ``proof,''
  ``reached the mathematical limit'') is maintained and enforced across
  all writeups.
\item \textbf{Adversarial review loop.}  Each iteration undergoes a
  structured adversarial review during development.
\item \textbf{Structured output.}  All numerical results are stored in
  machine-readable JSON with full-precision strings.  No result is
  reported from terminal output alone.
\end{itemize}

\end{document}